%

%
%
%

\input amstex
\documentstyle{amsppt}
\hsize = 5.4 truein
\vsize = 8.7 truein
\baselineskip = 12pt
\NoBlackBoxes
\TagsAsMath
\define\kint{(\kappa,\ \kappa^+)}

\define\si{\sigma}
\define\omu{\overline\mu}
\define\oA{\overline{A}}
\define\olp{\overline{p}}
\define\olb{\overline\beta}
\define\orht{\overset\rightharpoonup\to}
\define\ok{\overline\kappa}
\define\la{\lambda}
\define\k{\kappa}
\define\vp{\varPhi}
\define\a{\alpha}
\define\be{\beta}
\define\de{\delta}
\define\De{\Delta}
\define\ga{\gamma}

\define\tf{\tilde f}
\define\tC{\tilde C}

\define\s{{\mu^\prime,\ 0,\ s}}
\define\od{\overline\delta}
\define\os{\overline s}
\define\orr{\overline\rho}
\define\on{\overline\nu}
\define\ol{\overline\lambda}
\define\sm{\setminus}
\define\Si{\Sigma}
\define\CS{\Cal S}
\define\th{\theta}
\define\N{\Cal N}
\define\M{\Cal M}
\define\om{\omega}

\define\nar{\narrower\smallskip\noindent}

\define\modseq{(\Cal N_i|\ i \leq \theta)}

\define\pr{\prime}

%

\topmatter
\title The Combinatorics of Combinatorial Coding by a Real\\
\endtitle
\author Saharon Shelah$^{1,2,4,5}$ and Lee J. Stanley$^{3,4,5}$
\endauthor
\thanks {\roster
\item"{1.}" Research partially supported by NSF,
the Basic Research Fund, Israel Academy of Science and MSRI\item"{2.}"
Paper number 425\item"{3.}"
Research partially supported by NSF grant DMS 8806536, MSRI and
the Reidler Foundation.
\item"{4.}" We are very grateful to the referee for his rare
combination of diligence and patience,
as well as for many helpful observations.
\item"{5.}" We would also
like to thank many colleagues
for their interest
and the administration and staff of MSRI and
for their hospitality during 1989-90.
\endroster%
}
\endthanks

\address{Hebrew University, Rutgers University}\endaddress
\address{Lehigh University}\endaddress
\abstract{We lay the combinatorial foundations for \cite{5}
by setting up and proving the essential properties of the
coding apparatus for singular cardinals.  We also prove another
result concerning the coding apparatus for inaccessible cardinals.}
\endabstract
\endtopmatter

\subheading{\S 0. INTRODUCTION}
\bigskip
In this paper, we lay the combinatorial foundations for the work
of \cite{5}.  For the most part, this involves setting up the
coding apparatus for singular cardinals, and proving its essential
properties, most notably the result about the existence of supercoherent
sequences, Lemma 3 (the Lemma of (1.4) of \cite{5}).
The sole exception occurs in
(11.2), where, as promised in (2.1.1) of \cite{5},
we show that we can assume some additional properties for the
system of $b_\a$, with $card\ \a$ inaccessible.

The combinatorial apparatus for singular cardinals is based on our work
in Part I, where, working in $L$, we prove that the
\lq\lq {\it Squarer Scales}"
principle holds.  This is Theorem 1, below; the proof stretches across
\S\S 1 - 6.  This material is based on (and improves) that of
\cite{2}; \cite{2} bears the same relationship to the material of
\cite{1}, which is where many of the basic ideas of
this construction
made their first published appearance.
\S 2, in particular, reviews the constuctions of \S\S 1 - 2 of
\cite{2}, without proofs.
In \S 7, we quote a \lq\lq classical" result of Jensen, from \cite{3}
which, again in $L$, gives a square system on singular cardinals.  The last
section of Part I hints at things to come in that it steps outside of $L$
to remark that the methods of \S\S 1 - 6 allow us
to build \lq\lq local versions" of the combinatorial system of (1.2) between
$\tau$ and $\tau^{+\omega}$ working in $L[X_\tau]$, where $X_\tau \subseteq
\tau^+$.  As noted there, the $X_\tau$ we have in mind are the
$A \cap \tau^{+\omega}$, where $A$ is as given by Lemma 3 of \cite{5}.

In Part II, we assume that $V = L[A]$, for this $A$ (and that $0^\sharp$
does not exist).  We show, in \S 9, how to transfer the combinatorial systems
of Part I to $V$.  This culminates in (9.4), where we define a
{\bf fine system of squares and pseudo-scales}
to be one which satisfies properties (A) - (D); these are
restatements of similarly labelled items of (1.2) of \cite{5}.  We
observe (the crucial fact having already been noted in (9.3.2)) that the
system obtained in (9.1) - (9.3) is indeed a fine system.
This is Corollary 2, below.

Of course, it is here that we make essential use of the Covering Lemma.
This is used to guarantee that $L$ \lq\lq gets the successors of
singular cardinals (cardinals of $V$) right", but also that our
$L\text{-scales}$ remain something close enough to $V\text{-scales}$.
Ostensibly, what is required for Part II is that if $d$ is a club subset
of $\k$, a singular limit of limit cardinals, such that $o.t.\ d < inf\ d$,
and $g \in V$ is a function with domain $d$ such that for
$\la \in d,\ g(\la) < \la^+$, then there is a function $f$ in our
$L\text{-scale}$ such that for sufficiently large $\la \in d,\ g(\la) <
f(\la)$.  In fact, something more is needed for the result of (11.1), namely
that the preceding holds when $d$ is any Easton set.
While it is \lq\lq folklore" that this follows from
the Covering Lemma when $0^\sharp$ does not exist, it is tempting, but
false, to think that this remains true without restriction to an Easton set,
as the referee pointed out.  This theme of restriction to an Easton set
is also implicit in \S 10 (the restriction to \lq\lq controlled" cardinals,
see (10.1.2) and (10.3)).

In \S 10, we prove Lemma 3, below, the Lemma of (1.4) of
\cite{5}, which states, roughly, that for the system of \S 9, there are enough
supercoherent sequences.
This is the centerpiece of this paper, and, in many ways of \cite{5} as well,
as the whole approach to \cite{5}, the precise formulation of the
definition of the forcing conditions, for example, was driven by
the plan of using Lemma 3 to underly the proof of distributivity.

Lemma 3 is proved in two stages, first, by proving, in
(10.2), that there are enough strongly coherent sequences, and then, in
(10.3) - (10.5), that if $(\Cal N_i|i \leq \th)$ is
strongly coherent then $(\Cal N_{\om i}|\ \leq \th)$ is supercoherent.
The arguments of (10.3) - (10.5) use the
most intricate properties of the system of
\S\S 1 - 6.  In \S 11, we close by proving two other, smaller results needed
in \cite{5}: in (11.1), we prove the Proposition of (1.5) of \cite{5} which
plays an important role in the proof, in (4.3) of \cite{5} that the
\lq\lq very tidy" conditions are dense, and in (11.2) we prove the
result mentioned above about the $b_\a$ for $\a$ such that $card\ \a$ is
inaccessible.

Before stating Theorem 1, Corollary 2
and Lemma 3, we should say that our notation is intended to
either be standard or have a clear meaning, or is introduced as needed.
It may, however, be worth pointing out that we use the same definitions of
$U(\k)$, for $\k$ a limit cardinal, as in \cite{5}; for singular
$\k$ this is reintroduced in (9.2); for inaccessible $\k$ this is reintroduced
in (11.2).  One instance where notation is required to do double duty is
$S_\a$; throughout most of Part I, this is the notion introduced in (1.1), but
on two occasions in the proof of (3.1), explicitly noted when they occur,
the same notation refers to Jensen's auxiliary hierarchy of \cite{3}.

We turn now to the statements of our main results.

\proclaim{Theorem 1} ($V = L$)  The {\it Squarer Scales} principle of
(1.2), below, holds.
\endproclaim

\proclaim{Corollary 2} ($0^\sharp$ does not exist and $V = L[A]$,
where $A$ is as given by Lemma 3 of \cite{3})  There is fine system
of squares and pseudo-scales, i.e., one satisfying (A) - (D)
of (9.4).
\endproclaim

\proclaim{Lemma 3} ( ... as in Corollary 2 ... )  The system of Corollary 2
has the additional property that whenever $\Cal M,\ \nu$ and $\theta$ are as
in (10.1), below, and $C \subseteq [H_{\nu^+}]^\theta$ is club then there is
super $\Cal M\text{-coherent } (\Cal N_i|i \leq \theta)$, with each
$|\Cal N_i| \in C$.
\endproclaim

%

\subheading{PART I:  LIFE IN $L$}
\bigskip

In Part I, comprising \S\S 1 - 8,
we develop the $L\text{-combinatorics}$ summarized in the {\it Squarer
Scales} principle of \S 1.
This is a strengthening of the
Squared Scales principle from \cite{2}.
In \S 2 we review material from \S\S 1, 2 of
\cite{2}.  In \S 3 we pause to give a more explicit
(and perhaps clearer) development of certain crucial
ideas implicit in \S\S 2, 3 of
\cite{2}; we then return to reviewing the material of \S 3 of
\cite{2}.
In \S 4 we introduce
a new fine structure parameter, and prove some its important properties.
Finally in \S 5, we rework the construction of \S\S 2 - 3 of \cite{2}
based on this new parameter, and
we prove the important lemmas which are the analogues of those
of \S 4 of \cite{2}.  This culminates, in \S 6, in the proof of:
\bigskip
\proclaim{Theorem 1}  In $L$, {\it Squarer Scales} holds.
\endproclaim
\bigskip
In \S 7, we recall Jensen's construction from \cite{3} of
a square system defined on ordinals, which, in $L$, are singular cardinals.
Finally, in \S 8, we note that the techniques of \S\S 1 - 6 allow us to
construct \lq\lq local versions" of the squares and scales obtained there.
More precisely, if $\tau = \aleph_2$, or $\tau$ is a limit cardinal, if
$X_\tau \subseteq \tau^+$ is such that, letting $\mu = \tau^{+\om},\
H_\mu = L_\mu [X_\tau]$, then, in $L[X_\tau]$,
working as in \S\S 1 - 6, we construct a scale between $\mu$ and $\mu^+$,
and for cardinals, $\la$, with $\tau < \la \leq \mu$, a square system
between $\la$ and $\la^+$, which will enjoy all the properties of
the \lq\lq global" system constructed in $L$.  In fact, we will not really
need all of the properties, here, and notably, not the Condensation Coherence
properties, but the construction is the same, and they fall out anyway.
Of course, the $X_\tau$ we have in mind are the $A \cap \tau^{+\om}$, where
$A$ is as guaranteed by \cite{4}, and as in Lemma 3 of \cite{5}.

As in \cite{2}, it will simplify notation if
we assume, throughout \S\S 1 - 7, that $\bold{V = L}$.  As there, however,
this is purely a matter of notational convenience.

\subheading{\S1.  SQUARER SCALES}

We state {\it Squarer Scales},
and point out how it is stronger than the principle of \cite{2}.
We state the strengthened principle in a notation
designed to be suggestive of that of \cite{5}
rather than in the notation of \cite{2}.
Thus, we write $\tf_{\om\nu}$  where $\vp^\nu$ was used in \cite{2}, etc.
We have however, kept
the same organization of items as in (4.11), of \cite{2}.   The
principal difference in the principles is that our (B)(5) is stronger than
that of (4.11) of \cite{2},
as our (B)(5) handles the $g \in \Cal S (\k)$ (see below) and
not just subfunctions of such $g$ whose domains are cofinal subsets of $\k$ of
small cardinality.  We need some preliminary definitions, which carry over to
the rest of Part I.
\bigskip
\noindent
\proclaim{(1.1)\ \ Definition}
$S$ will denote the class of
ordinals, $\nu$, such that there is $\om < \a < \om\nu$ for which
$J_\nu \models$ \lq\lq $\a$ is the largest cardinal".
For $\nu \in S,\ \a_\nu$ is the unique such $\a$.
$S_\a$ is $\{\nu \in S|\a_\nu = \a\}$.  For limit cardinals,
$\ \k,\ \CS(\k)$ is the set of functions, $g$, such that
$dom\ g$ is a final segment of
the uncountable cardinals smaller than $\k$, and for
$\ok \in dom\ g,\ g(\ok) \in (\ok,\ \ok^+)$.
As usual, if $f,\ g \in \CS(\k),\ f <^* g$ if
for some $\k_0,\ \om < \k_0 < \k$ for all cardinals, $\ok$ with $\k_0 \le
\ok < \k,\ \tilde f(\ok) < \tilde g(\ok)$ and $f \leq^* g$ iff the final
$\lq\lq <$" is replaced by $\lq\lq \leq$".
\endproclaim

We should note that the above is the \lq\lq official" definition of
$S_\a$, but that in \S 3, below, we use this notation for a different notion.
This will noted when it occurs.
\bigskip
\noindent
{\bf (1.2)\ \ THE PRINCIPLE.}
\medskip
{\narrower\medskip\noindent
There is a sequence $(C_{\om\nu}|\nu \in S)$,
and for each limit cardinal, $\k$,
a sequence $(\tf_{\om\nu}|\nu \in S_\k$ \&\ o.t. $C_{\om\nu} <
\k)$ such that:\medskip}

\roster
\item"{(A)}"  For all $\nu\in S$, letting $\a = \a_\nu$:

\itemitem {(1)}  $C_{\om\nu}$ is a closed subset of
$\{\om\tau|\tau \in S_\a \} \cap \om\nu;\ sup\ C_{\om\nu} <
\om\nu \Rightarrow cf\ \om\nu = \om$,
\itemitem {(2)}  $\om\ga \in C_{\om\nu} \Rightarrow
C_{\om\ga} = \om\ga \cap C_{\om\nu}$,
\itemitem {(3)}  o.t. $C_{\om\nu} \le \a$, and if $\a$ is a singular cardinal,
then for sufficiently large $\nu \in S_\a, <$ holds.
\medskip

\item"{(B)}"  For all limit cardinals, $\k$, all
$\nu \in S_\k,\ o.t.\ C_{\om\nu} < \k \Rightarrow \tf_{\om\nu} \in \CS(\k)$ and:

\itemitem {(1)}  $\ok \in dom\ \tf_{\om\nu} \Rightarrow
(\tf_{\om\nu}(\ok)$ is a limit ordinal and
$\on \in S_{\ok}$,
where $\om\on = \tf_{\om\nu}(\ok))$,
\itemitem {(2)}  $\ok \in dom\ \tf_{\om\nu} \Rightarrow
(\forall\om\la\in C_{\om\nu})\ok \in
dom\ \tf_{\om\la}$, and $\tf_{\om\la}(\ok) \leq \tf_{\om\nu}(\ok)$,
\itemitem {(3)}  $(\forall\tau \in S_\k \cap \nu) o.t.\ C_{\om\tau} < \k
\Rightarrow \tf_{\om\tau} <^* \tf_{\om\nu}$,
\itemitem {(4)}  Suppose that $sup\ C_{\om\nu} = \om\nu \ \&
\ \ok \in dom\ \tf_{\om\nu}$.  If
$(\tf_{\om\la}(\ok)|\om\la \in C_{\om\nu})$ is not eventually constant
then $\tf_{\om\nu}(\ok) = sup\ \{\tf_{\om\la}(\ok)|\om\la \in C_{\om\nu}\}$,
\itemitem {(5)}  if $\k$ is singular, then
whenever $g \in \Cal S(\k)$, there is
$\nu_0 \in S_\k$ such that $o.t.\ C_{\om\nu_0} < \k$ and $g <^* \tf_{\om\nu_0}$.
\medskip

\item"{(C)}"  For limit cardinals, $\k,\ \& \ \nu \in S_\k$, if
$o.t.\ C_{\om\nu} < \k$ and $\ok \in
dom\ \tf_{\om\nu}$, then, letting $\om\overline\nu = \tf_{\om\nu}(\ok)$ and
$\Phi = \{\tf_{\om\la}(\ok)|\om\la \in C_{\om\nu}\}$:

\itemitem {(1)}  $\Phi$ is a final segment of
$C_{\om\overline\nu}$ (we take this to include the case where $C_{\om\on}$ is
bounded in $\om\on$ and $\Phi = \emptyset$),
\itemitem {(2)}  $\tf_{\om\overline\nu} = \tf_{\om\nu}|\ok$,
\itemitem {(3)}  $\Phi \in J_{\overline\be}$, whenever
$J_{\overline\be} \models$ \lq\lq $\om\overline\nu$ is not a cardinal".
\endroster
\medskip
\proclaim{Remark}  We only use the scales for $\k$ which are singular
cardinals, but the construction gives them for inaccessibles as well.
In \S 9, we ignore the scales for inaccessibles.
\endproclaim
\bigskip
\subheading{\S 2.  REVIEW OF \S\S 1 - 2 OF \cite{2}}
\bigskip
\noindent
{\bf(2.1) THE COLLAPSING STRUCTURES.}
\medskip
\proclaim{(2.1.1)  Definition}
For $\nu \in S$, if $\om\nu$ is not a cardinal,
$\beta(\nu)$ is the least $\beta \ge \nu$ such that
$J_{\be + 1} \models$ \lq\lq $\om\nu$ is not a cardinal".
\endproclaim

Let  $\be = \be(\nu)$; then, for some $n$
there is $f$, which is $\Si_{n+1}$-definable over $J_\be$ (in parameters from
$J_\be$) and $f$ is a map onto $\om\nu$ from a
subset of a smaller ordinal.
\medskip
\proclaim{(2.1.2)  Definition}
$n(\nu)$ is the least $n$ such that there is such an $f$ which is
$\Si_{n+1}$-definable over $J_\be$ (in parameters from $J_\be$).
Let $n = n(\nu)$; then $\rho(\nu)$ is $\rho^n_\be$, the $n^{th}$-projectum of
$\be,\ A(\nu) = A^n_\be =$ the $n^{th}$-master-code of $\be$,
and setting $\rho =
\rho(\nu),\ A = A(\nu),\ \Cal A(\nu) = (J_\rho,\ \in,\ A)$.
It can be shown that
$\rho^{n+1}_\be \le \a_\nu$ and $\nu \le \rho$, so that for some finite set of
ordinals, $p \subseteq \om\rho$, all elements of $J_\rho$ are $\Si_1$-
definable in $\Cal A(\nu)$ (i.e., are unique solutions in $\Cal A(\nu)$ of
$\Si_1$-formulas in one free variable) using parameters from $\a_\nu \cup p$.

We abbreviate this last assertion by writing: $J_\rho =
h\lq\lq (\om \times (\a_\nu \cup p))$,
where $h = h_{\Cal A(\nu)} = h_\nu$
is the canonical $\Si_1$-Skolem function for $\Cal A(\nu)$.  We let
$p(\nu)$ be the least such $p$ with respect to the lexicographic ordering
of the \underbar{decreasing} enumeration of $p$.  Then, $\Cal A^+(\nu) = (\Cal
A(\nu),\ p(\nu))\ (p(\nu)$ is a new individual constant).
This is the collapsing structure for $\nu$.
\endproclaim
\medskip
An important and useful
fact is provided by Corollary (1.8) of \cite{2}:
if $n = 0$ then there is a
largest cardinal $\ga$ in the sense of $J_\be$, and $p(\nu) \not\subseteq \ga$;
further, if $(X,\ \in) \prec_{\Si_1} (J_\be,\ \in)$ and $p(\nu) \in X$, then
$\ga \in X$.
\bigskip
\noindent
{\bf(2.2)\ \ CLOSING THE CLASS OF COLLAPSING STRUCTURES.}
\medskip
We close off the class of collapsing structures under
transitive collapses of (constructible) rudimentarilly-closed
substructures.
\medskip
\proclaim{Definition}
$\Cal O^+ := \{\Cal A^+(\nu)|\nu \in S,\ \om\nu \text{ is
not a cardinal} \}$.  We let $(\Cal B,\ q) \in
\tilde{\Cal O^+}$
iff $|\Cal B|$ is transitive and for some $(\Cal A,\ p) \in \Cal
O^+,\ (\Cal B,\ q)$ is isomorphic to a (constructible) $rud(\Cal A)$-closed
substructure of $(\Cal A,\ p)$.
\endproclaim

Thus, if
$(\Cal B,\ q) \in \tilde{\Cal O}^+,\ \Cal B$ is
amenable and of the form
$(J_{\rho^\prime},\ \in,\ A^\prime)$.  Further,
$\tilde{\Cal O}^+$ is closed for taking transitive collapses of constructible
$rud(\Cal B)$-closed substructures of
$(\Cal B,\ q) \in \tilde{\Cal O}^+$; in particular, it is
closed under amenable initial segments and transitive collapses of
constructible $\Si_1$-elementary substructures.
\bigskip
\noindent
{\bf (2.3) A SQUARE SYSTEM ON} $\bold{\tilde{\Cal O}^+}$.
\medskip
For $s = (\Cal A,\ p) = (J_\rho,\ \in,\ A,\ p) \in \tilde{\Cal O}^+$,
a closed subset,
$\tC_s \subseteq \om\rho$ is constructed; $\tC_s$ is cofinal in $\om\rho$ if
$cf\ \om\rho > \om$.  Crucial in the definition and structure of
$\tC_s$ are
the sets $\De(\xi,\ s)$ for $\xi < \om\rho$, where $\om\de \in \Delta(\xi,\ s)$
iff $\om\de < \om\rho$ and for some $\be,\ \om\de = sup\ h_s\lq\lq (\om \times
(\be \cup \{\xi\}))$.  Recall that for
a set $X$ of ordinals,\ $X^\prime$ is the set of limit points of $X$, below
$sup\ X$.  First, consider $\De(0,\ s)$: if this is empty,
$cf\ \om\rho = \om$ and
$\tC_s = \emptyset$.  If this is cofinal in $\om\rho$, then $\tC_s =
(\De(0,\ s))^\prime$; of course, if $\De(0,\ s)$ is cofinal in $\om\rho$, then
$\tC_s$ is cofinal in $\om\rho$ if
$cf\ \om\rho > \om$.
The remaining case is when $\De(0,\ s)$ has a largest element, $\om\de$.  Then,
for some $\be,\ \om\de = sup(OR \cap h_s\lq\lq (\om \times \be))$, but
$\om\rho = sup(OR \cap h_s\lq\lq (\om \times (\be + 1)))$; note that
this can occur even if $cf\ \om\rho > \om$, since we must consider all the
unique solutions in $s$ of $\Si_1$ formulas
$\phi(\nu_0,\ \xi_1,\ \cdots,\ \xi_n,\ \be)$, where
$\xi_1,\ \cdots,\ \xi_n < \be$; so all we have
for certain is that $cf\ \om\rho \leq cf\ \be$.

In this case, we set
$\be = \be^0_s,\ \de = \de^1_s(\de^0_s = 0$, for all $s$).
We have the same trichotomy for $\De(\om\de^1_s,\ s)$: if
$\De(\om\de^1_s,\ s) =
\emptyset,\ cf\ \om\rho = \om$ and in this case,
$\tC_s = (\De(0,\ s))^\prime$; if
$\De(\om\de^1_s,\ s)$ is cofinal in $\om\rho$ then $\tC_s = (\De(0,\ s))^\prime
\cup (\De(\om\de^1_s,\ s))^\prime$.
Finally, if $\De(\om\de^1_s,\ s)$ has a largest element, $\om\de$,
then we have $\be = \be^1_s,\ \de = \de^2_s$ such that $\om\de^2_s$ is the
largest element of $\De(\om\de^1_s,\ s),\ \om\de^2_s =
sup(OR \cap h_s\lq\lq (\om \times (\be \cup \{\om\de^1_s\})))$ and $\om\rho =
sup(OR \cap h_s\lq\lq (\om \times ((\be + 1) \cup \{\om\de^1_s\})))$.

The crucial observation, proved in (2.40) of \cite{2},
is that, in this case, $\be^0_s >
\be^1_s$.  Thus, the process terminates after a finite number, $m_s \ge 1$, of
steps; in all cases, $\tC_s =
\bigcup\{(\De(\om\de^i_s,\ s))^\prime|i < m_s\}$.
If $i = m_{s-1},\ \de = \de^i_s$, then $\tC_s$ has a (possibly empty) final
segment, $(\De(\om\de,\ s))^\prime$ and if $cf\ \om\rho > \om$ then
$(\De(\om\de,\ s))^\prime$ is cofinal in $\om\rho$, since otherwise
$\De(\om\de,\ s)$ would have a largest element (the other possibilities are
eliminated by the cofinality hypothesis), which is impossible since the process
terminates after $m_s$ steps.

It is not really necessary to \lq\lq thin out"
by taking only the limit points of the $\De(\om\de^i_s,\ s)$, but this slightly
facilitates the proof of the coherence property of the $\tC_s:$ if $\om\de \in
\tC_s$, then, setting $s^\prime = s|J_\de,\ \tC_{s^\prime} =
\tC_s \cap \om\de$.  As an important preliminary step it is shown that if, for
$t \in \tilde {\Cal O}^+$, we let $a_t = \{\om\de^i_t|i < m_t \}$, then, for
all $s \in \tilde {\Cal O}^+$ and all $\om\de \in \tC_s$, letting $s^\prime
= s|J_\de,\ a_{s^\prime} = a_s \cap \om\de$.
Of course, if we chose not to thin out, then the coherence property would
hold for $\om\de \in (\tC_s)^\prime$, and we could, by choosing constructible
cofinal $\om$-sequences in the appropriate cases, guarantee that $\tC_s$ is
\underbar{always} cofinal in $\om\rho$.  Jensen has taken this approach in
\cite{1}, where the cofinal $\om$-sequences
are chosen in a canonical and natural fashion.
\bigskip
\noindent
\subheading{\S 3.  A CLOSER LOOK AT THE $X_{\s}$, AND \S 3 OF \cite{2}}
\medskip
We prove three Lemmas related to the structure of the $X_{\s}$.  The first, in
(3.1), guarantees that when
$\mu$ is a singular cardinal, $\nu \in S_\mu$ and $\nu$ is
sufficiently large that $\Cal A(\nu) \models \lq\lq \mu$ is singular", then,
letting $s = \Cal A^+(\nu)$, for some
$\mu^\prime < \mu,\ X_\s$ is cofinal in $\om\rho$.  This is certainly
well known
to fine-structure experts, but was never stated explicitly in
\cite{2}.  For completeness, we give it here.
Some of the ideas involved in (3.1) and (3.2)
appear in the proofs of (4.1) and (4.3) of \cite{2}.

The second Lemma, in (3.2) shows that when $\mu^\prime$ is as guaranteed by
(3.1), then, under two additional, mild assumptions, $\tilde C_s \subseteq
X_\s$.  The third Lemma, (3.3), explores what occurs when $X_\s$ is not cofinal
in $\om\rho$.  Essentially, it shows that if $s^\prime =
s|\de^{\prime\prime} \in \tilde{\Cal O}^+$, then, at least as far as
$X_{s^\prime,0,\mu^\prime}$ is concerned, we can assume without loss of
generality that either $\de^{\prime\prime} = \rho(s)$ or that
$\de^{\prime\prime} \in X_\s$.  These Lemmas will be heavily used in
\S\S 5, 6, below.

\proclaim{(3.1)\ \ Lemma}  Assume that
$\mu$ is a singular cardinal, $\nu \in S_\mu$ and $\nu$ is
sufficiently large that $\Cal A(\nu) \models \lq\lq \mu$ is singular",
and let $s = \Cal A^+(\nu)$.  Then,
for some $\mu^\prime < \mu,\ X_{\s}\
\text{is cofinal in}\ \om\rho(\nu).$
\endproclaim
\demo{Proof} Let $f: a \to_{onto}\ \om\nu$ be $\Si_1(\Cal A(\nu))$ in
parameters $\orht y \in J_{\rho(\nu)}$.  Suppose, e.g., that $\phi$ is a
$\Si_1$ formula such that $\zeta = f(\xi) \Leftrightarrow \Cal A(\nu) \models
\phi(\zeta,\ \xi,\ \orht y)$.  Let $\theta$ be the $\Si_0$ formula such that
$\phi$ is $\exists v_0\theta(v_0,\ \zeta,\ \xi,\ \orht y)$.  Let
$\theta^\prime(\eta,\ \xi)$ be:

$$\orht y \in S_\eta \land (\exists v_0 \in S_\eta)(\exists \zeta <
\eta)\theta(v_0,\ \zeta,\ \xi,\ \orht y).$$

\noindent
{\bf In the above,} $\bold{S_\eta}$ {\bf is the}
$\bold{\eta^{\text {th}}}$ {\bf stage in
Jensen's auxiliary hierarchy, not the notion defined in (1.1), above}.
Note that if $\theta^\prime(\eta,\ \xi)$ and $\eta < \eta^*$, then
$\theta^\prime(\eta^*\ \xi)$.
Let $g(\xi) \simeq\ $ the least $\eta$ such that $\theta^\prime(\eta,\ \xi)$.
Thus, $g$ is $\Si_0$.

Well known arguments (involving the downward extension of embeddings Lemma)
then show:

$$(\forall \mu^\prime \leq \mu)[\om\nu = sup\ f\lq\lq (a \cap \mu^\prime)
\Leftrightarrow \om\rho = sup\ g\lq\lq (a \cap \mu^\prime)].$$

\noindent
Thus, if there is $\mu^\prime < \mu$ with $\om\rho =
sup\ g\lq\lq (a \cap \mu^\prime)$, (3.1) holds, so, towards a contradiction,
assume that $g\lq\lq (a \cap \mu^\prime)$ is bounded in $\mu$ for all
$\mu^\prime < \mu$.  For such $\mu^\prime$, let $\si(\mu^\prime) =
sup\ g\lq\lq (a \cap \mu^\prime)$.  Also, let $g_{\mu^\prime} = g|\mu^\prime$,
so $g_{\mu^\prime} \subseteq \mu^\prime \times \si(\mu^\prime)$ and
$g_{\mu^\prime}$ is definable over $S_{\si(\mu^\prime)}$ (Jensen's auxiliary
hierarchy again).

This makes it clear that each $g_{\mu^\prime} \in
J_{\rho(\nu)}$, and, in fact that $\mu^\prime \mapsto g_{\mu^\prime}$ is
$\Si_1(\Cal A(\nu))$ in parameters $\orht y$.  But then, the same holds for
$\mu^\prime \mapsto \si(\mu^\prime)$, and,
denoting this last function by $\si,\ \si$ is non-decreasing with domain $\mu$.
Now, let $g^* \in J_{\rho(\nu)}$ be a map of a subset of some $\mu^\prime <
\mu$ cofinally into $\mu$.  Then, the function $\si \circ g^*$ is
$\Si_1(\Cal A(\nu))$ and maps a subset of $\mu^\prime$ cofinally into
$\om\rho(\nu)$, contradiction.  This completes the proof.
\enddemo

\proclaim{Remark}
If $\mu,\ \nu,\ s$ are as in (3.1), then (3.1) clearly gives that
$o.t.\ \De(0,\ s) < \mu$, and therefore, for all $\xi \in a_s,
\ o.t.\ \De(\xi,\ s) < \mu$.  But then, clearly $o.t.\ \tC_s < \mu$.
\endproclaim

\proclaim{(3.2)\ \ Lemma}
If $\mu,\ \nu$,
etc., are as in (3.1), then whenever $\mu^\prime < \mu$ is as guaranteed by
(3.1), $(a_s \subseteq X_{\s}\ \&\ o.t.\ \tC_s \leq \mu^\prime)
\Rightarrow \tC_s \subseteq X_{\s}$.
\endproclaim

\demo{Proof}
Let $f_{\s} =
(\os,\ |f|,\ s)$, where $|f|$ is the inverse of the transitive collapsing
map for $s|X_{\s}$.  We first argue in the case where
$\tC_s$ is cofinal in $\om\rho$.  Then, applying
(2.31)(b) of \cite{2}
to $f_{\s}$, we get that $\tC_{\os}$ is cofinal in
$\om\orr$, where $\orr = \rho(\os)$.  But then, since $range\ |f|$ is
cofinal in $\om\rho$, in fact, $|f|\lq\lq \tC_{\os}$ is cofinal in
$\om\rho$, and by (2.31)(c) of \cite{2},
$|f|\lq\lq \tC_{\os} \subseteq
\tC_s$, so $X_{\s} \cap \tC_s$ is cofinal in $\om\rho$.  But then, let
$\om\de \in X_{\s} \cap \tC_s$.  Since $\tC_{s|\de}$ is an initial segment of
$\tC_s,\ o.t.\ \tC_{s|\de} < \mu^\prime$.  Finally, by (2.25) of \cite{2}
(whose statement contains a typo; the statement should read:
$\lq\lq \cdots,\ \de < \rho(s)$ and $s|\de \in \tilde {\Cal O}^+$ then
$\cdots$ "), it easily follows that $\tC_{s|\de} \in X_{\s}$.  But then, since
$o.t.\ \tC_{s|\de} < \mu^\prime$, in fact $\tC_{s|\de} = \tC_s \cap \om\de
\subseteq X_{\s}$.  Thus, arbitrarily large initial segments of $\tC_s$ are
included in $X_{\s}$.

This completes the proof when $\tC_s$ is
cofinal in $\om\rho$.  When $\tC_s = \emptyset$, there is nothing to prove.
So, suppose that $\tC_s$ has a greatest element.  Since $X_{\s}$ is cofinal in
$\om\rho$, it follows from (2.31) and (2.38) of \cite{2}
that if $f_{\s}$ is as above,
then $\tC_{\os}$ has a largest element and that $|f|(max\ \tC_{\os}) =
max\ \tC_s$.  Then, arguing as above, and appealing, once again, to
(2.31) and also (2.38) of \cite{2}, the conclusion is clear.
\enddemo
\bigskip
\noindent
{\bf (3.3)}
\medskip
In dealing with the situation where
$X_{\s}$ is not cofinal in $\om\rho(s)$,
it will facilitate some
of the arguments to replace $\mu^\prime$ by $\mu^\prime + 1$, so that,
letting $f = |f_{\mu^\prime + 1,\ 0,\ s}|,\
f(\mu^\prime,\ \mu^\prime + 1) =
(\mu^\prime,\ \mu^\prime + 1)$.  This also is faithful to the context in which
we shall apply this material, in \S\S 5, 6, below.
We adopt the same notation as in (3.1) and (3.2),
but with $\de = \de(f_{\mu^\prime + 1,\ 0,\ s}) < \rho(s)$.
\bigskip
\proclaim{Lemma}  Whenever $\de^{\prime\prime} \leq \rho(s)$ and
$s|\de^{\prime\prime} \in \tilde {\Cal O}^+$, there is $\de^* \in
\{\rho(s)\} \cup X_{\mu^\prime + 1,\ 0,\ s}$ such that
$s|\de^* \in \tilde {\Cal O}^+$ and
$|f_{\mu^\prime + 1,\ 0,\ s|\de^{\prime\prime}}| =
|f_{\mu^\prime + 1,\ 0,\ s|\de^*}|$.
\endproclaim

\demo{Proof}
Suppose, first, that $\de \leq \de^{\prime\prime} <
\rho(s)$.  By (2.23) of \cite{2}, $f =
|f_{\mu^\prime + 1,\ 0,\ s|\de^{\prime\prime}}|$.  Next,
suppose that $f(\om\od^*) = \om\de^* > sup\ f\lq\lq \om\od^*$, where
$\os|\od^* \in \tilde {\Cal O}^+$.
Let $g = (\os|\od^*,\ f|J_{\od^*},\ s|\de^*)$.  By (2.32) of \cite{2},
$s|\de^* \in \tilde {\Cal O}^+$ and
$f|J_{\od^*}:\os|\od^*\ \to_{\Si_\om}\
s|\de^*$.
But then, we
clearly have that $X_{\mu^\prime + 1,\ 0,\ \os|\od^*}$ is cofinal in
$\om\od^*$ iff $X_{\mu^\prime + 1,\ 0,\ s|\de^*}$ is cofinal in $\om\de^*$.
However, since $X_{\mu^\prime + 1,\ 0,\ s|\de^*} \subseteq
X_{\mu^\prime + 1,\ 0,\ s} \cap
J_{\de^*}$, and $sup\ f\lq\lq \om\od^* =
sup\ (X_{\mu^\prime + 1,\ 0,\ s} \cap J_{\de^*})$,
clearly $X_{\mu^\prime + 1,\ 0,\ \os|\od^*}$ is not cofinal in $\om\od^*$,
i.e., $\de(f_{\mu^\prime + 1,\ 0,\ \os|\od^*}) < \od^*$.
Let $\od^\prime = \de(f_{\mu^\prime + 1,\ 0,\ \os|\od^*})$, and let
$\de^\prime = |f|(\od^\prime)$.  By (2.30) of \cite{2}, $\de^\prime =
\de(f_{\mu^\prime + 1,\ 0,\ s|\de^*})$.

Then, $\de^\prime < \de(g) = sup\ f\lq\lq \om\od^*$.  Finally,
if $\de^\prime \leq \de^{\prime\prime} < \de^*$, applying
(2.23) of \cite{2}, with $s|\de^*$ in place of $s,\ g$ in place of $f,
\ \de^\prime = \de(g)$ in place
of $\de(f)$ (in the notation of (2.23) of \cite{2}, $\de = \de(f)$), we have
$|f_{\mu^\prime + 1,\ 0,\ s|\de^{\prime\prime}}| =
|f_{\mu^\prime + 1,\ 0,\ s|\de^*}|$.  This completes the proof.
\enddemo
\bigskip
\noindent
{\bf (3.4)\ \ \lq\lq PROJECTING" A TAIL OF } $\bold{\tilde C_s }$ {\bf TO A
SUBSET OF } $\bold{\om\nu.}$
\medskip
In \S 3 of \cite{2},
$C_{\om\nu}$ is defined, for $\nu \in S,\ \om\nu$ not a
cardinal.  First, a final segment of $\tC_s$ is chosen, where $s = \Cal
A^+(\nu)$.

\proclaim {(3.4.1) Definition}  Let $\om\de \in\tC_s,\ \om\de > \a_\nu$,
let $s^\prime = s|J_\de$, and
let $Y = Y_{\de,\ \nu} = h_{s^\prime}\lq\lq (\om \times \a_\nu)$; then,
$\om\de \in \hat C_\nu$ iff $\a_\nu \in Y$.
\endproclaim

It is shown in (3.2)(b) and (3.3) of \cite{2}
that if $\om\rho = \om\nu$, then $\hat C_\nu = \tC_s$.
For $\om\de \in \hat C_\nu$, it is shown, in (3.2)-(3.4) of \cite{2},
that there is unique $\la$ such that $\la\in S_\a$ and $\Cal A^+(\la) \cong
s^\prime|Y_{\de,\ \nu} = s|Y_{\de,\ \nu}$.

\proclaim{(3.4.2) Definition}  For $\om\de \in \hat C_\nu$ we set
$\la(\de,\ \nu) = $ the unique $\la \in S_\a$
such that $\Cal A^+(\la) \cong
s^\prime|Y_{\de,\ \nu} = s|Y_{\de,\ \nu}$.
\endproclaim

If $\rho = \nu$, then $\la = \de$, as is shown in (3.3) of \cite{2}.
An important observation is made in (3.2)(a) of \cite{2}:
$Y_{\de,\ \nu}$ is cofinal in $\om\de$.

\proclaim{(3.4.3) Definition} $C_{\om\nu} =
\{\om\la(\de,\ \nu)|\om\de \in \hat C_\nu\}$.
\endproclaim

It is then shown in (3.6) - (3.8) that the $C_\nu$ have the correct
properties, i.e., those of (A) of (1.2), above.
\bigskip

\subheading{\S 4.  A NEW PARAMETER}

Our main tool in proving the strengthened version, (B)(5) of (1.1), above,
of the (B)(5) of \cite{5}, is a small
but potentially quite useful Lemma, below, involving
a new parameter which we now introduce.
Then, in \S 5, we supply the arguments which
replace those of \S 4 of \cite{2},
making the changes and improvements enabled by this Lemma.

\proclaim{(4.1)\ \ Definition}
Let $\nu \in S,\ \rho = \rho(\nu),\ \Cal A = \Cal A(\nu),
\ \a = \a_\nu,\ \a \le \tau \le \om\rho$.  Let $R_\nu(p,\ \tau)$ be the
property: $p \in [\om\rho]^{<\om} \&
\ h_{\Cal A}\lq\lq (\om\times(\tau \cup p)) = J_\rho$.
Let $P_\nu(p)$ be the
property $R_\nu(p,\ \a)$; let $Q_\nu(p)$ be the property
$R_\nu(p,\ \a + 1)$.  So
$p(\nu) = $ the least $p$ such that $P_\nu(p)$,
with respect to lexicographic order of the
decreasing enumeration of finite subsets of $\om\rho$.  Analogously,
\underbar{define} $q(\nu) =$ the least $q$ such that $Q_\nu(q)$,
with respect to the same ordering.
\endproclaim

\proclaim{Remarks}

\roster
\item $p(\nu) \cap \a_\nu = \emptyset$; $q(\nu) \cap \a_\nu + 1 = \emptyset$,
\item $q(\nu) = p(\nu) \Leftrightarrow \a_\nu \in
h_{\Cal A(\nu)}\lq\lq (\om \times(\a_\nu\cup q(\nu)))$,
\item $\forall r(P_\nu(r) \Rightarrow Q_\nu(r))$, so $q(\nu) \le_* p(\nu)$,
\item $P(\{\a_\nu\} \cup q(\nu));$ thus $p(\nu) \le_* \{\a_\nu\}\cup q_\nu$.
\endroster
\endproclaim

\proclaim{(4.2)\ \ Lemma}  Either $p(\nu) = q(\nu)$ or $p(\nu) = q(\nu)
\cup \{\a_\nu\}$.
\endproclaim

\demo{Proof} Let $p = p(\nu),\ q = q(\nu),\ \a = \a_\nu$.
Note that by Remark 3, if $p \neq q$, then
there is $c \subset p$ which is a common final segment of
$p$ and $q$ and either $c = q$ or else the largest member
of $q \setminus c$ is less than the largest member of $p \setminus
c$.  However, by Remark 4, there is $d \subseteq q \cup \{\a\}$ which
is a common final segment of $p$ and $q \cup \{\a\}$, and if
$d \neq p$, then $d \neq q \cup \{\a\}$ and the largest member
of $p \setminus d$ is less than the largest member of $(q
\cup \{\a\}) \setminus d$.  In the latter case
the largest member of $(q \cup \{\a\}) \setminus d$ must be greater
than $\a$ so it is simply the largest member of $q \setminus d$, and we have a
contradiction.  Thus, we must have that $d = p$.  If $\a \not\in p$,
then $p \subseteq q$, which is also impossible.  Thus, $\a \in p$ and
$p = q \cup \a$.
\enddemo

The main difference between the arguments in \S\S 5, 6, below,
and those of \S 4 of \cite{2} is that for $s = \Cal A^+(\nu)$, below, we
use $X^*_{\mu^\prime,\ 0,\ s} = X_{\mu^\prime + 1,\ 0,\ s}
\ = h_s\lq\lq (\om \times (\mu^\prime +1))$,
whereas, in \S 4 of \cite{2}, we used $X_{\mu^\prime,\ 0,\ s} =
h_s\lq\lq (\om \times \mu^\prime)$.  Of course,
$X_{\s} \subseteq X^*_{\s}\ \&\ \mu^\prime \in X^*_{\s}$.
By the above Lemma, either $p(\nu) \neq q(\nu)$, in which case, we have
$X_{\s} = X^*_{\s}$, or else $\a_\nu \not\in p(\nu)$.  One main observation is
that none of this really depends on $\mu^\prime$.
\bigskip
\subheading{\S 5.  REWORKING \S 4 OF \cite{2}}
\medskip
In this Section, we rework the material corresponding to (4.5) - (4.10)
of \cite{2}.  There is no analogue of (4.8), however, because of our use
of the $X^*_{\s}$.  (5.1) corresponds to (4.5) of \cite{2}.  (5.2) corresponds
to (4.7), of \cite{2}, in ideas, if not in statement.  In (5.3), we define
the $\tilde f_{\om\nu}$ (the analogous definitions in \cite{2} were
(4.6) and (4.9)).  (5.4) corresponds to (4.10) of \cite{2} and establishes the
Condensation Coherence property, (C)(1) of (1.2).

\proclaim{(5.1) Lemma}  If $\nu \in S_\mu,\ \mu$ is a limit cardinal, $s =
\Cal A^+(\nu)$, if $\mu^\prime < \mu,\ \mu^\prime$ is a cardinal and
$f^*_{\s} = (\os,\ |f|,\ s)$, where $|f|:
\os \to s|X^*_{\mu^\prime,\ 0,\ s}$ is the inverse of the transitive
collapsing map, then there is a unique $\on \in S_{\mu^\prime}$, such that
either $\om\on = OR \cap |\os|$ or
$\os \models$ \lq\lq $\om\on$ is a cardinal";
further, either $X_{\mu^\prime,\ 0,\ s} =
X^*_{\mu^\prime,\ 0,\ s}\ \&
\ \os = \Cal A^+(\on)$, or $X_{\mu^\prime,\ 0,\ s} \neq
X^*_{\mu^\prime,\ 0,\ s},\ \os = (\Cal A(\on),\ q(\on))\ \&
\ \mu^\prime \in p(\on) \sm q(\on)$.
\endproclaim

\demo{Proof}
The existence and uniqueness of $\on$ are immediate from the fact that
$\mu^\prime < sup\ X^*_{\s}$.  To get the remainder of the Lemma, we shall
apply the downward extension of embeddings lemma to $|f|$.  Let $n = n(\nu)$,
and let $\os = (J_{\orr},\ \oA,\ \olp)$.
The downward extension of embeddings  gives us a $\olb$ and
$\hat f:\ J_{\olb} \to_{\Si_{n + 1}} J_{\be(\nu)},
\ |f| \subseteq \hat f$, such that $\orr = \rho^n_{\olb}$ and
$\oA = A^n_{\olb}$.  Since
$h_{\os}\lq\lq (\om \times (\mu^\prime + 1)) = J_{\orr}$, as usual
we have $\olb = \be(\on)$ and $n \geq n(\on)$.

For the reverse inequality, if $n = 0$, there is nothing to prove, so
suppose $n > 0$.  Then, if
$n > n(\on)$, exactly as in (3.1), (3.3) and (3.4) of \cite{2},
we would have $\rho(\on) \geq \rho^{n(\on) + 1}_{\olb}
\geq \orr$, on the one hand, but $\rho^{n(\on) + 1}_{\olb}
\leq \mu^\prime + 1 < \on \leq \orr$, on the other hand, a
contradiction.

Thus, $\os = (\Cal A(\on),\ \olp)$.  Of
course, $|f|(\olp) = p(\nu)$ and, by construction, $J_{\orr} =
h_{\Cal A(\on)}\lq\lq (\om \times ((\mu^\prime + 1) \cup \olp))$;
i.e. $Q_{\on}(\olp)$ holds.  If $q \in J_{\orr},\ q <_*
\olp$, and for some $i < \om$ and
$\overset\rightharpoonup\to\xi \in
[(\mu^\prime + 1)]^{<\om},\ \olp =
h_{\Cal A(\on)}(i,\ \overset\rightharpoonup\to\xi,\ q)$, then $p =
h_{\Cal A(\nu)}
(i,\ |f|(\overset\rightharpoonup\to\xi),\ |f|(q))$ and $|f|(q) <_* p$.
This, however, contradicts the definition of $p = p(\nu)$ since
$|f|(\overset\rightharpoonup\to\xi) =
\overset\rightharpoonup\to\xi \in [\mu]^{<\om}$ (recall that
$|f|(\mu^\prime) = \mu^\prime)$.
Thus, $\olp = q(\on)$.  By (4.2), either $\olp = p(\on)$ and $J_{\orr} =
h_{\Cal A(\on)}\lq\lq (\om\times((\mu^\prime \cup \olp))$, in which
case $\mu^\prime \in X_{\s}$, so $X_{\s} = X^*_{\s}$ and $\os = \Cal
A^+(\on)$, or $p(\on) \neq \olp$, in which case $\mu^\prime
= \a_{\on} \notin \olp,\ p(\on) = \{\a_{\on}\} \cup \olp =
\{\mu^\prime\} \cup \olp$.  Then, $\mu^\prime \notin
h_{\Cal A(\on)}\lq\lq (\om\times(\mu^\prime \cup \olp))$, so
$\mu^\prime \in X^*_{\s} \sm X_{\s},\ \os = (\Cal A(\on),
\ q(\on))$.  Note, here, that $\mu^\prime \cup p(\on) =
(\mu^\prime + 1) \cup \olp = (\mu^\prime + 1) \cup q(\on)$.
\enddemo
\bigskip


\proclaim{(5.2) Propostion}  Let $\nu,\ \mu,\ s$ be as in (5.1).
Assume that $\hat C_\nu \neq \emptyset$, and let $\om\de_0 = inf\ \hat C_\nu$.
In addition to our hypotheses on $\mu^\prime$ from
(5.1), suppose further that $\mu \in X_{\mu^\prime + 1,\ 0,\ s|\de_0}$.
Let $X^*_{\s},
\ f^*_{\s},\ |f|,\ \os$, etc. be as in (5.1).  Let
$\on$ be as guaranteed by (5.1).
Let $\om\de \in \hat C_\nu$.  Let $Y = Y_{\de,\ \nu}$ and let
$\la = \la(\de,\ \nu)$.  Suppose that $\de \in X^*_{\s}$, and let
$|f|(\od) = \de$.  Let $\ol$ be as guaranteed by (5.1) with $\la$ in the place
of $\nu$.  Then, $\ol = \la(\od,\ \on)$.
\endproclaim

\demo{Proof}  Note that our additional hypothesis on $\mu^\prime$ guarantees
that the analogous statement holds for {\it any} $\om\eta \in \hat C_\nu$ in
place of $\om\de_0$.  The only real difficulty
in proving the Lemma is that, in general, $X^*_{\mu^\prime,\ 0,\ s|\de}
\subset X^*_{\mu^\prime,\ 0,\ s} \cap J_\de$.

Let $\pi: (J_{\rho^\prime},\ A^\prime,\ p^\prime) \to
(s|\de)|Y$ be the isomorphism,
and let $s^\prime = (J_{\rho^\prime},\ A^\prime,\ p^\prime)$
so, by \S 3 of \cite{2},
$s^\prime = \Cal A^+(\la)$.  As remarked after (3.4.2),
$\pi\lq\lq \om\rho^\prime$ is cofinal in $\om\de$.
By (5.1), $\overline s = (\Cal A(\overline\nu),\ q(\overline\nu))$,
with the dichotomy of the conclusion of (5.1).

Clearly,
$X^*_{\mu^\prime,\ 0,\ \os|{\od}} = |f|^{-1}[X^*_{\mu^\prime,\ 0,\ s|\de}]$ and
so, letting $\overline
Y = X^*_{\mu^\prime,\ 0,\ \os|{\od}},\ \overline Y$ is cofinal in $\om\od$;
this follows immediately from (2.30) and (3.2)(a) of \cite{2}.
Also, here we have $\mu^\prime \in \overline Y$.

The following easy
observation will be important in establishing (B)(4) and (C)(1) of (1.2),
above; this will be done in (5.4), below:

{\nar ($\ast$) $\ \ \overline Y = X_{\mu^\prime,\ 0,\ \overline t}$, where
$\overline t = (\Cal A(\on),\ p(\on))|J_{\od}$.  So, depending on the
dichotomy of (5.1), either $\overline t = \os|J_{\od}$ or
$\overline t = (\Cal A(\on),\ \{\mu^\prime\} \cup \olp)|J_{\od}$.
\medskip}

\noindent
It is clear from ($\ast$) that $\om\od \in \hat C_{\on}$, in either case.

Let $\overline\pi:
(J_{\orr^\prime},\ \overline A^\prime,\ \overline p^\prime) \to
(\os|\od)|\overline Y$ be the inverse of the transitive collapse.  So
$\os^\prime\ :=\
(J_{\orr^\prime},\ \overline A^\prime,\ \overline p^\prime) =
(\Cal A(\ol),\ q(\ol))$,
and either $p(\ol) = q(\ol)$ or $p(\ol) = \{\mu^\prime\} \cup
q(\ol)$.  In either case, $(\os^\prime,\ |f| \circ \overline\pi,\ s|\de) =
f^*_{\mu^\prime,\ 0,\ s|\de}$.
Then, clearly, $range\ |f| \circ \overline \pi \subseteq Y$, and
$(\os^\prime,\ \pi^{-1} \circ |f| \circ \overline\pi,\ \os) =
f^*_{\mu^\prime,\ 0,\ \os}$.  It then follows easily that $\ol =
\la(\od,\ \on)$, as required.
\enddemo
\bigskip

\proclaim{(5.3) Definition}  Let $\mu,\ \nu,\ s$ be as in (5.1).
Let $\mu^*_1(\nu) =\ $ the least uncountable cardinal,
$\mu^\prime < \mu$, such that for all $\om\de \in \hat C_\nu,
\ \mu \in X^*_{\mu^\prime + 1,\ 0,\ s|\de}$.  Thus, if $\hat C_\nu = \emptyset,
\ \mu^*_1(\nu) = \aleph_1$.  Otherwise,
as remarked at the beginning of the proof
of (5.2), this is just the least
$\mu^\prime$ such that $\mu \in X^*_{\mu^\prime,0, s|\de_0}$, where
$\de_0$ is as in (5.2).  For cardinals $\mu^\prime
\in [\mu^*_1(\nu),\ \mu)$, let
$\tilde f_{\om\nu}(\mu^\prime) = \om\on$, where $\on$ is as guaranteed by
(5.1).
\endproclaim

\proclaim{(5.4) Proposition} (Condensation Coherence):
If $\nu,\ \mu,\ s$, etc., are as
in (5.1) and $\mu^*_1(\nu) \le \mu^\prime < \mu,\ \mu^\prime$ a cardinal,
then letting $\om\overline\nu = \tilde f_{\om\nu}(\mu^\prime)$:
\roster
\item "{(a)}" for all $\om\la \in C_{\om\nu},
\ \mu^*_1(\la) \le \mu^\prime$,
\item "{(b)}" let $\Phi =
\{\tilde f_{\om\la}(\mu^\prime)|\om\la \in C_{\om\nu}\}$; then $\Phi$
is a final segment of $C_{\om\overline\nu}$
(we take this to include the case where $C_{\om\on}$ is bounded in $\om\on$ and
$\Phi = \emptyset$),
\item "{(c)}" if $\mu^\prime$ is a limit cardinal, then $\mu^\prime >
\mu^*_1(\nu),
\ \mu^*_1(\on) = \mu^*_1(\nu)$ and $\tilde f_{\om\overline\nu} = \tilde
f_{\om\nu}|\mu^\prime$.
\endroster
\endproclaim

\demo{Proof}  (a) is clear and (c) follows easily from ($\ast$) of (5.2), above.
For (b), let $s = \Cal A^+(\nu) = (J_\rho,\ A,\ p)$, and let
$f = f^*_{\s} = (\os,\ |f|,\ s)$.  Let $\de = \de(f)$.
If $range\ |f|$ is not cofinal in $\om\rho$, then as we have already noted in
arguing for (3.3), above, $X^*_{\s} =
X^*_{\s|\de} \subseteq Y_{\de,\ \nu}$,
so composing with $\pi^{-1}$, the isomorphism
between $s|\de$ and $\Cal A^+(\la(\de,\ \nu))$, we transport the whole
situation down to $\la^* = \la(\de,\ \nu)$.  Now, if (b) holds between
$\la^*$ and $\on$, as we shall argue that it will,
we can use (3.3), above, to
conclude that it holds between $\nu$ and $\on$, since (3.3) gives
that $\Phi = \{\tilde f_{\om\la}(\mu^\prime)|\om\la \in C_{\om\la^*}\}$.

Thus, we may assume that $\de = \rho$, i.e. that $range\ |f| = X^*_{\s}$ is
cofinal in $\om\rho$.  This allows us to appeal to (2.31) of \cite{2}
to conclude that,
letting $\os = (J_{\orr},\ \oA,\ \olp),\ \tC_s$ is cofinal in $\om\rho$ iff
$\tC_{\os}$ is cofinal in $\om\orr$, that $\tC_s = \emptyset$ iff $\tC_{\os} =
\emptyset$, and that if $\om\de,\ \om\od$ are the maxima of $\tC_s,
\ \tC_{\os}$, respectively, then $\de = |f|(\od)$.  Now, since $\mu^*_1(\nu)
\leq \mu^\prime,\ \mu \in X^*_{\s}$, so let $\omu$ be such that $\mu =
|f|(\omu)$.  Recalling the last clause of (5.1), above, it is then easy to see
that:

{\nar ($\ast$)\ \ $\hat C_\nu \neq \emptyset$
iff there is $\om\od \in \tC_{\os}$
such that $\omu \in
h_{\Cal A^+(\on)|J_{\od}}\lq\lq (\om \times \mu^\prime)$.
\medskip}

\noindent
Thus, if $\hat C_\nu = \emptyset$, then $C_{\om\on}$ is bounded in $\om\on$ and
$\Phi$ is the empty final segment of
$C_{\om\on}$.  So, for the remainder of the proof, we assume that
$\hat C_\nu \neq \emptyset$.

Let $\om\de \in \hat C_\nu$ and let
$\la = \la(\de,\ \nu)$.  By (3.3), above, we may suppose that, as in
(5.2), $\de = |f|(\od)$.
Adopt the notation of (5.2), above.  We proved there
that $\om\od \in \hat C_{\on}$ and that
$\tf_{\om\la}(\mu^\prime) = \la(\od,\ \on)$,
so $\tf_{\om\la}(\mu^\prime) \in C_{\om\on}$, for all
$\om\la \in C_{\om\nu}$.  It remains
only to show that letting $W =
\{\om\od|\om\de \in \hat C_\nu\ \cap X^*_{\s}\}$, then $W$ is a final
segment of $\hat C_{\on}$.  This, however, is clear, since $W =$\newline
$\hat C_{\on} \setminus \om\od_0$, where $|f|(\om\od_0) =
inf\ \hat C_\nu \cap X^*_{\s}$.
\enddemo
\bigskip

\proclaim{(5.5)  Remark}
We should point out that
$\tf_{\om\nu}(\mu^\prime) = \be(f^*_{\mu^\prime,\ 0,\ s})$.
\endproclaim
\bigskip
\subheading{\S 6.  COMPLETING THE PROOF OF SQUARER SCALES}
\medskip

(A) of (1.2) is immediate from the material of \S\S 1 - 3
of \cite{2}, summarized in \S\S 2, 3, above.
(B)(1) is clear from construction.  (B)(2) follows easily
from the definition of $\mu^*_1$, in (5.3), above,
the remark about $\mu^*_1$ in (5.3), above, prior to the definition of
$\tilde f_{\om\nu}$
and the proof of (3.3), above.
(C)(1) and (C)(2) follow easily from (5.4).  (C)(3)
follows from (5.4) and the analogous statement about $C_{\om\on}$,
but the latter follows readily from (2.25) and (2.33) of \cite{2}.

We argue for (B)(4).  Let $\om\on = \tf_{\om\nu}(\ok)$.
We should note that the hypothesis that $\Phi$ has limit
order type will hold if $X^*_{\ok,\ 0,\ s}$ is cofinal in $\om\rho(\nu)$, by
(3.2), above, where $s = \Cal A^+(\nu)$.  Let
$f^*_{\ok,\ 0,\ s} = (\os,\ |f|,\ s)$.  As in (5.1),
$\os = (\Cal A(\on),\ q(\on))$.
Applying (C)(1), we have that $\Phi$ is a final
segment of $C_{\om\on}$.  However, since $\Phi$ has limit order type, by
hypothesis, it must therefore be cofinal in $\om\on$.

It remains to verify the scale
properties, (B)(3) and (B)(5).
We first argue for (B)(5); we shall appeal to a part of its proof in arguing
for (B)(3).  So, let $\k$ be singular and let $g \in \Cal S(\k)$.  Clearly
there is $\nu_0 \in S_\k$ such that $g \in J_{\nu_0}$, and of course, taking
$\nu_0$ sufficiently large, we may suppose that
$J_{\nu_0} \models\ \lq\lq \k$\ is singular".
But then, as in the arguments for (3.1) and (3.2), above,
$o.t.\ C_{\om\nu_0} < \k$.  Since $J_{\nu_0}
\subseteq J_{\rho(\nu_0)}$, it will suffice to prove:

{\nar ($\ast$)  if $\k$ is a singular cardinal,
$\eta \in S_\k,\ o.t.\ C_{\om\eta} < \k,\ g \in \Cal S(\k)
\cap J_{\rho(\eta)}$, then $g <^* \tf_\eta$.
\medskip}

\demo{Proof of ($\ast$)}  Let $s = \Cal A^+(\eta)$
and let $\ok < \k$ be such that $g \in h_s\lq\lq (\om \times \ok)$.  Let
$\ok,\ \mu^*_1(\eta) \leq \k^\prime < \k$ be a cardinal.  We shall argue that
$g(\k^\prime) < \tf_\eta(\k^\prime)$.  The main observation is that since
$\k^\prime \in X^*_{\k^\prime,\ 0,\ s}$, we also have $g(\k^\prime) \in
X^*_{\k^\prime,\ 0,\ s}$.  But then, since $s \models \lq\lq card
\ g(\k^\prime) = \k^\prime$", clearly $g(\k^\prime) + 1 \subseteq
X^*_{\k^\prime,\ 0,\ s}$ and so $|f||(g(\k^\prime) + 1) =
id|(g(\k^\prime) + 1)$.  Thus,
letting $f^*_{\k^\prime,\ 0,\ s} =
(\os,\ |f|,\ s),\ \os \models \lq\lq card\ g(\k^\prime) = \k^\prime$" and so
$g(\k^\prime) < \be(f^*_{\k^\prime,\ 0,\ s}) = \tf_\eta(\k^\prime)$.  The last
equality is by (5.5), above.  This completes the proof of
($\ast$) and therefore of (B)(5).

We finish by arguing for (B)(3).  In view of ($\ast$),
and since $J_\nu \subseteq J_{\rho(\nu)}$, it will clearly suffice to
show that if $\tau \in S_\k \cap \nu,
\ o.t.\ C_{\om\tau} < \k$ then $\tf_{\om\tau} \in
J_\nu$.  Now, under these
hypotheses, it is clear that $\be(\tau) < \nu$, and therefore that
$\Cal A^+(\tau) \in J_\nu$ and so, letting
$s = \Cal A^+(\tau),\ h_s \in J_\nu$.  But
then, the function $\k^\prime \mapsto \be(f^*_{\k^\prime,\ 0,\ s})$ is also an
element of $J_\nu$.  Finally, in virtue of (5.5),
$\tf_{\om\tau}$ is a restriction of this function to the set of cardinals in
a final segment of its domain and therefore
$\tf_{\om\tau} \in J_\nu$, as required.
\enddemo
\bigskip
\noindent
\subheading{\S7.  A Square on Singular L-Cardinals}
\medskip
We simply recall that in \cite{3}, Jensen contructed a system
$(\tilde D_\k|\k\text{ is a singular } L\text{-cardinal})$ with the properties
that $\tilde D_\k \subset \k$ is a club of $\k$ such that
$o.t. \tilde D_\k < min D_k$ and such that
if $\la \in (\tilde D_\k)^\prime$, then $\la$ is a singular $L\text{-cardinal}$
and $D_\la = D_\k \cap \la$.
\bigskip
\noindent
\subheading{\S8.  Local Systems in $\bold{L[X_\tau]}$}
\medskip
Prior to (1.1), we outlined the thrust of this section, so we limit ourselves
to the statement of the result.
\proclaim{Lemma}  Suppose that $\tau = \aleph_2$ or $\tau$ is a limit
cardinal,
and let $\mu = \tau^{+\om}$.  Suppose that $X_\tau
\subseteq \tau^+$ is such that
$H_\mu = L_\mu [X_\tau]$.  Let
$S^\tau = \bigcup\{S_\la|\tau < \la \leq \mu\text{
and } \la\text{ is a cardinal }\}$.
Then, in $L[X_\tau]$, there are systems $(C_{\om\nu}|\nu \in S^\tau)$ and
$(\tilde f_{\om\nu}|\nu \in S_\mu\ \&\ o.t.\ C_{\om\nu} < \mu)$
which satisfy (A) - (C) of (1.2),except that, in addition,
we require that if $\la \in dom\ \tilde f_{\om\nu}$, then $\la > \tau$.
\endproclaim

Of course, the $X_\tau$ we have in mind are the $A \cap \tau^{+\om}$.
\bigskip\bigskip

%

\subheading{Part II:  Life In A Sharpless V}
\bigskip
In Part II, which comprises \S\S 9 - 11, we transfer the combinatorial
structures of Part I to a sharpless $V$, and prove the results required for
\cite{5}, notably Corollary 2 and Lemma 3 (Lemma (1.4) of \cite{5}).
As in \cite{5}, we work in the context provided by \cite{4}, i.e., we
assume that $0^\sharp$ does not exist and we
work in $L[A]$, where
$A \subseteq OR$ is such that for all uncountable cardinals $\kappa,
\ H_\kappa = L_\kappa[A]$, such that $A =  (A \cap \omega_2) \cup \bigcup\{A
\cap (\kappa, \kappa^+)|\kappa \in \Lambda\}$, where $\Lambda$
is the class of limit cardinals together with $\aleph_2$.
Further, if $\kappa = \aleph_2$, or $\kappa$ is inaccessible,
then for $\delta \in (\kappa,\ \kappa^+),\ (card\ \delta)^{L[A \cap \delta]} =
\kappa$.
\medskip
In \S 9, we show how to transfer the combinatorial systems of Part I to
$V$, indicating briefly how the necessary modifications are performed.
We culminate, in (9.4), with the definition of a {\bf fine system of squares
and pseudo-scales} and the observation that the system obtained in (9.1) -
(9.3) is indeed a fine system.  This proves Corollary 2 and
corresponds to (1.2) of \cite{5}.  In \S 10 we prove
Lemma 3.  We finish, in \S 11, by proving two smaller results, used in
(1.5) of \cite{5} and  (2.1.1) of \cite{5}.

{\bf In the remainder of this paper, notions such as \lq\lq cardinal",
\lq\lq singular cardinal", etc., mean \lq\lq cardinal in the sense of} $\bold
V\text{{\bf ",}}$ {\bf \lq\lq singular cardinal in the sense of $\bold
V\text{{\bf", etc.}}$
\bigskip
\noindent
\subheading{\S9.  From L to V}
\bigskip
\noindent
{\bf (9.1) OBTAINING THE } $\bold {D_\k}$ {\bf FROM THE } $\bold {\tilde
D_\k.}$
\medskip
First, for singular cardinals, $\mu$ of the form $\eta^{+\om}$, we let
$\Lambda$ be as above, we let
$\eta^*$ be the unique member of $\Lambda$ such that $\mu = (\eta^*)^{+\om}$,
and we define $D_\mu\ :=\ \{\aleph_\tau \in (\eta^*,\mu)|\tau\text{ is odd}\}$.

So, assume that $\k$ is a singular limit of limit cardinals.
Let $E_\k$ be the set of singular cardinals in $\tilde D_\k$, where
$\tilde D_\k$ is as in \S 7.  If
$(E_\k)^\prime$ is cofinal in $\kappa$,
let $D^*_\kappa =$\newline
$\kappa \cap (E_\k)^\prime$ and
set $\lambda \in I(\kappa)$ iff $\lambda$ is a successor point of
$D^*_\kappa$.  If $(E_\k)^\prime$ is
bounded in $\kappa$, set $D^*_\kappa = \emptyset,\ I(\kappa) = \{\kappa\}$.
Note that if $\lambda \leq \kappa,\ \lambda \in I(\kappa)$, then $\lambda$ is a
singular limit of limit cardinals and $\lambda \in I(\lambda)$.  Also, note
that if $\lambda \in I(\lambda)$, then $cf\ \lambda = \omega$.  Thus, for all
singular limits of limit cardinals, $\lambda$, such that $\lambda \in
I(\lambda)$, choose $x(\lambda) = \{\lambda_j|j < \omega \}$, cofinal in
$\lambda,\ (\lambda_j|j < \omega)$ increasing, such that:
\medskip
\roster
\item  min ${\tilde D}_\lambda < \lambda_0$; whenever
$\lambda^\prime \in \lambda \cap E^\prime,
\ \lambda^\prime < \lambda_0$,
\item  for all $j < \omega$, there is $\delta(j)$, which is
not a successor ordinal, such that $\lambda_j =
\aleph_{\omega(\delta(j)+j+1)}$.
\endroster
\bigskip
Then, for all $\kappa$ which are singular
limits of limit cardinals, let $D_\kappa = D^*_\kappa \cup
\bigcup\{x(\lambda)|\lambda \in I(\kappa)\}$.
\bigskip
Note that by construction, $(D_\kappa|\kappa$ a singular limit of
limit cardinals$)$ has the usual coherence property; further, letting
$\delta_\kappa = o.t.\ D_\kappa,\ \delta_\kappa \leq o.t. \tilde D_\kappa <
min\ {\tilde D}_\kappa < min\ D_\kappa$, and letting $(\lambda^\kappa_i|i <
\delta_\kappa)$ increasingly enumerate $D_\kappa$, if
$\lambda^{\kappa_1}_{i_1} = \lambda^{\kappa_2}_{i_2}$, then $i_1 = i_2$.
Further, note that for all $\lambda,\ \{\delta_\kappa|\lambda \in D_\kappa\}
\subseteq \lambda$.
\bigskip
\noindent
{\bf (9.2)  MODIFYING THE } $\bold {C_\alpha}.$
\medskip
If $\k$ is a singular cardinal, then, by Covering, $\k^+ = (\k^+)^L$, so
that the system $(C_{\om\nu}|\nu \in S_\k)$ is very close to being a
square-system between $\k$ and $\k^+$.
In fact, in virtue of (3.1) and (3.2), above, except for an initial segment,
$I$, of $\alpha \in S_\kappa$,
we always have $o.t.\ C_\alpha < \kappa$.
Recall that
as in \cite{5}, for singular $\k$, we let $U(\k)$ be the set of multiples of
$\k^2$ in $(\k,\ \k^+)$.
Let $\phi_\k$ be the continuous
order-isomorphism between
$\{\om\nu|\nu \in S_\k\ \setminus I\}$
and the set of limit multiples of $\k^2$ in
$(\k,\ \k^+)$.  We transfer the system to live on the latter set, via $\phi$, by
taking $C_{\phi(\a)}\ :=\ \phi\lq\lq C_\a$.
Finally, the $C_\alpha$ constructed in Part I
are not necessarily club:  they have been
thinned by removing successor points.  These are restored, in a canonical
way by recursion on the well-founded relation \lq\lq $\a \in C_\be$" by
supplying cofinal $\om\text{-sequences}$ above $sup\ C_\a$ to those
$\a$ whose $C_\a$ is not cofinal.  We have abused notation
by using $C_\a$ to denote
this modified system as well.
\bigskip
\noindent
{\bf (9.3)  MODIFYING THE } $\bold{\tilde f_{\om\nu}.}$
\medskip
There are several kinds of modifications we carry out.  The first is to
transfer the scales to live on the $(U(\k))^\prime$, as we did for the
squares, in (9.2).  Here, it is a bit more complicated, since we must also
transfer the values, via different continuous order-isomorphisms.  Also,
at least in the first few stages of the modifications, we continue to deal
with certain $L\text{-cardinals}$ which may not be cardinals of $V$.

So, if $\k$ is an $L\text{-cardinal}$, we let $\phi_\k$ be the
order-isomorphism of
$\{\a \in S_\k|o.t.\ C_\a < \k\}$ to an initial segment,
$T_\k$, of the set of
limit multiples of $\k^2$.  Note that if $\k$ is actually
a cardinal, then $T_\k = (U(\k))^\prime$. Further, if $\k$ is
actually a singular cardinal, then $\phi_\k$ is as
in (9.2).  Finally, if $\k$ is actually a regular cardinal,
then $\phi_\k$ is only $<\ \k\text{-continuous}$ but, as
will be clear, that is all that is required.

\proclaim{(9.3.1)  Definition}
Now, suppose that $\k$ is actually a singular limit of limit cardinals.
We define
$\hat f_\eta$ for $\eta \in (U(\k))^\prime$, with domain the set of
$L\text{-cardinals}$ between $\aleph_1$ and $\k$.  Let $\a \in S_\k$
with $o.t.\ C_\a < \k$
be such that $\eta = \phi_\k(\a)$.  First, suppose that
$\la \in dom\ \tilde f_{\om\a}$.  We then set $\hat f_\eta(\la)\ :=\
\phi_\la(\tilde f_{\om\a}(\la))$.  If $\aleph_1 < \la < \k,\ \la$
is an $L\text{-cardinal}$ and $\la \not\in dom\ \tilde f_{\om\a}$,
we set $\hat f_\eta(\la)\ :=\ \la^2\om$.

If $\tau = \aleph_2$ or $\tau$ is a limit cardinal and
$\mu = \tau^{+\om}$, the procedure is similar:  for $\eta \in (U(\mu))^\prime$,
letting $\a \in S_\mu$ with $o.t.\ C_\a < \mu$ be such that $\eta =
\phi_\mu(\a)$, if $\la \in dom\ \tilde f_{\om\a}$, we set
$\hat f_\eta(\la)\ :=\ \phi_\la(\tilde f_{\om\a}(\la)$, but we only extend the
domain to be the set of cardinals between $\tau$ and $\mu$, again, using
$\la^2\om$ as the default value.

Next, we must define the scale functions $\hat f_\eta$, for
$\eta \in U(\k) \setminus (U(\k))^\prime$, where $\k$ is a singular cardinal.
This is rather straightforward.
First, if $\k$ is a singular limit of limit cardinals,
let $\la \in X\ iff \la$ is an $L\text{-cardinal}$ and
$\aleph_1 < \la < \k$,
so suppose that $\k$ is an $\om\text{-successor}$.
If $\k = \aleph_\om$, let $\tau = \aleph_2$; otherwise,
let $\tau$ be the unique limit cardinal with $\k = \tau^{+\om}$.
In both of these cases, let $\la \in X\ iff \tau < \la < \k$
and $\la$ is a cardinal.
if $0 < n < \om$ and $\eta= \k^2n$,
for all $\la \in X$,
we let $\hat f_\eta(\la)\ :=\ \la^22n$.  Otherwise, let $\si$ be a limit
ordinal, $0 < n < \om$, and suppose that $\eta = \k^2(\si + n)$.  Then, for
all $\la \in X$, we set
$\hat f_\eta(\la)\ :=\ \hat f_\si(\la) + \la^22n$.
\endproclaim

\proclaim{(9.3.2)  Remark}  It is easy to see that the transferred system
of $C_\eta$ and $\hat f_\eta$ for
$\eta \in (U(card\ \eta))^\prime$ satisfies
the obvious analogues of (A) - (C) of (1.2), above.  We shall use
this observation in
(9.4) and in \S 10, without additional comment.
\endproclaim

\proclaim{(9.3.3)  Definition}
Finally, we define the $f^*_\eta$ for $\eta \in U(\k)$, where $\k$ is a
singular cardinal.  These are simply $\hat f_\eta|D_\k$, where $D_\k$ is as
given by (9.1).
\endproclaim
\bigskip
\noindent
{\bf (9.4)  A FINE SYSTEM.}
\medskip
We now define the notion of a {\bf fine system of squares and pseudo-scales}
as one which satisfies properties (A) - (D), below (these are restatements
of the similarly labelled items of (1.2) of \cite{5}).  When this is done,
it will be clear (by (9.3.2)) that
since we are assuming that $0^\sharp$ does not exist and that
$V = L[A]$, where $A$ is as given by Lemma 3
of \cite{5}, the combinatorial system developped in (9.1) - (9.3)
is a fine system of squares and pseudo-scales.
This proves Corollary 2.
\medskip
\proclaim{Definition}  A {\bf fine system of squares and pseudo-scales} is
a system $(D_\mu|\mu$ is a singular limit of limit cardinals),
$(C_\alpha|\alpha \in (U(\k))^\prime \cap \k^+ \& \k\ \text{is a singular
cardinal})$, $(f^*_\alpha|\alpha \in U(\k),\ \& \k\ \text{is a singular
cardinal})$ satisfying the following properties (A) - (D).
\endproclaim

\noindent
\roster
\item"{(A)}"
\medskip
\noindent
For singular cardinals, $\mu,\ D_\mu$ is a club subset of the set of
cardinals less than
$\mu$ such that if $\mu$ is a limit of limit cardinals, then
all members of $D_\mu$ are singular, while if
(($\tau = \aleph_2$ or $\tau$ is a limit cardinal) and
$\mu = \tau^{+\om})$, then $\la \in D_\mu\ iff (\tau < \la < \mu\
\& \la = \aleph_\xi$, where $\xi$ is odd), and:
\medskip
\itemitem {(1)}  $o.t.\ D_\mu < min\ D_\mu$,
\itemitem {(2)}  if $\lambda$ is a limit point of $D_\mu,\
D_\lambda = D_\mu \cap \lambda$.
\itemitem {(3)}  if $\lambda \in D_\mu$ is not a limit point
of $D_\mu$ then $\lambda$ is not a limit of limit cardinals.
\itemitem {(4)}  suppose that $\lambda \in D_{\k_i},\ i = 1, 2$, and let
$j_i$ be such that $\lambda$ is the $j_i^{th}$ member of $D_{\k_i}$.  Then,
$j_1 = j_2$.
\medskip
\item"{(B)}"  For singular cardinals, $\k$,
and $\alpha \in (U(\k))^\prime \cap \k^+),\ C_\alpha$ is a club
subset of the set of even multiples of $\k^2$ below $\a$,
of order type less than $\k$, and such that if
$\be \in C_\a$ but is not a limit point of $C_\a$, then $\be$ is not
a limit point of $U(\k)$, and with the usual coherence property:
if $\be$ is a limit point of $C_\alpha,\ C_\beta = C_\alpha \cap \beta$.
\medskip
\item"{(C)}"  For singular cardinals, $\k$, and
$\alpha \in U(\k)),\ dom\ f^*_\alpha = D_\k$,
for $\la \in D_\k,\ f^*_\a(\la)$ is an even multiple of $\la^2$
and:
\itemitem {(1)}  if $\kappa < \alpha < \beta,\ \a,\ \be \in U(\k)$
then $f^*_\alpha <^* f^*_\beta$, i.e., for some $\lambda_0 < \kappa$,
whenever $\lambda \in D_\kappa \setminus \la_0,\
f^*_\alpha(\lambda) < f^*_\beta(\lambda)$; further,
if $\alpha \in C_\beta$, then the preceding holds
for {\bf all } $\lambda \in D_\kappa$,
\itemitem {(2)}  whenever $g$ is a function with $dom\ g = D_\k$
and for all $\la \in D_\k,\ g(\la) < \la^+$, for some $\a \in U(\k),\
g <^* f^*_\a$,
\itemitem {(3)}  if $\k$ is a singular limit of limit cardinals,
$\la \in D_\k,\ \a \in U(\k),\ \a^\prime = f^*_\a(\la)$
and $\la^\prime \in D_\k \cap \la$, then $f^*_\a(\la^\prime)
= f^*_{\a^\prime}(\la^\prime)$, and if $\k$ is not a limit of
limit cardinals and $\a,\ \be \in U(\k),\ \la \in D_\k$ and
$f^*_\a(\la) = f^*_\be(\la)$, then $f^*_\a|\la = f^*_\be|\la$,
\itemitem {(4)}  for limit points, $\a$, of $U(\k)$, and
$\la \in D_\k$, $\Phi(\a,\la)\ :=\ \{f^*_\be(\la)|\be \in C_\a\}$
is a final segment of $C_{f^*_\a(\la)}$; further,
on a tail of $D_\k,\ \Phi(\a,\ \la)$ has limit order type.
\endroster
\medskip
We recall the observation made in (1.2) of \cite{5} to the effect that
even though the $f^*_\a$ are not defined when $card\ \a$
is a successor cardinal,
nevertheless the property of the second clause of (3) allows us to define them
in a conventional way so that we will then have the property of the first
clause of (3), even for $\k$ which are not limits of limit cardinals.
\bigskip
\roster
\item"{(D)}"  Decodability of (A) - (C):  For all singular $\k,\ D_\k$
and the systems $(C_\a|\a < \k^+\text{ is a limit point of } U(\k)),\
(f^*_\a|\a \in U(\k))$ are
canonically definable in $L[A \cap \k]$.
\endroster
\medskip
To make it completely clear why this
follows from (9.3.2), it will be useful
to give the correspondence between items of (B), (C), above, and the
items of the {\it Squarer Scales} principle of (1.2).  (B) corresponds to (A)
of (1.2).  (C)(1) corresponds to the conjunction of (B)(2) and
(B)(3) of (1.2).  (C)(2) corresponds to
(B)(5) of (1.2).  (C)(3) corresponds to (C)(2) of (1.2).  (C)(4) corresponds
to the conjunction of (B)(4) and (C)(1)of (1.2).  (D) corresponds to
(C)(3) of (1.2).
\bigskip

%

\subheading{\S 10.  THE EXISTENCE OF SUPER-COHERENT SEQUENCES}
\bigskip
In this section we prove Lemma 3 (Lemma (1.4) of \cite{5}).
This lemma states
that there are \lq\lq enough" super-coherent sequences.
We do this by first showing, in
(10.2), that there \lq\lq enough" strongly coherent sequences, and then, in
(10.5), showing that if $(\Cal N_i|i \leq \th)$
is strongly coherent then $(\Cal N_{\om i}|i \leq \th)$  is super-coherent.
The proofs of (10.3) - (10.5) exploit the most subtle
combinatorial properties of the {\it Squarer Scales}.  For convenience, we
begin by restating the definitions of strongly coherent and super-coherent,
and some preliminary related notions from (1.1) and (1.3) of \cite{5}.
Following (10.2) we lay out the plan for the proof carried out in (10.3)
- (10.5).
\bigskip
\noindent
{\bf (10.1)\ \ MODEL SEQUENCES AND COHERENCE.}
\medskip
Let $\th > \aleph_1$ be regular.  Let $\Cal M = (H_{\nu^+},\ \in,\ \cdots, )$,
where $\nu$ is a singular cardinal, $\nu >> \th$ and $(H_\nu,\ \in)$ models a
sufficiently rich fragment of ZFC.
Let $\si \leq \th$ and let
$(\Cal N_i:i \leq \si)$ be an increasing
continuous elementary tower of elementary substructures of $\Cal M$.

\proclaim{(10.1.1)  Definition}
We say that $(\Cal N_i|i \leq \si)$ is
$\bold{(\Cal M,\ \th)\text{{\bf-standard }}}$ {\bf of length }
$\bold{\si + 1}$
if, letting $N_i \ :=\  |\Cal N_i|$, for all $i \leq \si,
\ card\ N_i = \th,\ \th + 1 \subseteq N_0$, for $i < \si,
\ [N_{i + 1}]^{<\ \th} \subseteq N_{i + 1}$ and, if $i$ is even,
$\Cal N_i \in N_{i + 1}$.
\endproclaim

\proclaim{(10.1.2)  Definition}
For such $\theta > \aleph_1$ and $\Cal M = (H_{\nu^+},\in,\cdots)$,
suppose that $\Cal N \prec \Cal M$, where,
letting $N \ :=\  |\Cal N|,
\ card\ N = \th$, and let $\k$ be a cardinal with $\th \leq \k,\ \k \in
N$.  Let $\chi_{\Cal N}(\k) = sup(N \cap \kint)$.

Recall that an Easton set of ordinals is one which is bounded below any
inaccessible cardinal.
For such $\Cal N$ and singular cardinals, $\k$, with
$\theta < \kappa \leq \nu$,
we say that $\kappa$ is $\bold {\Cal N-\text{{\bf{controlled}}}}$ if
there is an Easton set $d$ with $\kappa \in d \in N$.

We define $p\chi_{\Cal N}$, an analogue of $\chi_{\Cal N}$, defined
on all singular cardinals, $\kappa$, which are $\Cal N-\text{controlled}$.
The definition makes sense for all cardinals $\k \in [\th,\ \nu]$, but we will
only use it for the singulars which are $\Cal N-\text{controlled}$.
If $\kappa \in N$, then of course $\kappa$ is $\Cal N-\text{controlled}$ and
in this case,
$p\chi_{\Cal N}(\kappa)\ :=\ \chi_{\Cal N}(\kappa)$.
Otherwise, $p\chi_{\Cal N}(\kappa)\ :=\ sup\ (\kappa^+ \cap Sk_{\Cal
M}(\{\kappa\} \cup N))$.
\endproclaim

The reason that we only consider controlled $\k$ is that (10.3), below,
gives an alternative characterization of
$p\chi_\N(\k)$ which is central in proving (10.5).
The alternative
characterization is equivalent only for controlled $\k$.
As we noted in \cite{5},
the restriction to such $\kappa$ is benign, for our
purposes.

\lq\lq Characteristic" functions of a model $\Cal N$
like $\chi_\Cal N$ and $p\chi_\Cal N$ often appear in
the work of the first author in a slightly different formulation, defined
to be \lq\lq pressing down" functions:  the value at a cardinal $\k$
is the supremum {\bf below} $\bold \k$ of some set of ordinals associated
with $\Cal N$.  Thus, in this formulation, our $\chi_\Cal N (\k)$ and
$p\chi_\Cal N (\k)$ would become values at $\k^+$ of these functions, and
we would also have at our dispostion the corresponding suprema below limit
cardinals.  In this connection, see the second Remark, follwoing the proof of
the Proposition in (10.3).

\proclaim{(10.1.3)  Definition}
Suppose $\k \in |\M|,\ \k$ is a singular cardinal, $\N,\ N$ are as
in (10.1.2), and
$\k \not\in N$.  Let $\mu_N(\k) =$ the least ordinal,
$\xi \in N$, such that $\xi > \k$ (clearly such exists, since
$\nu \in N$).  Clearly $\mu_N(\k)$ is a limit of limit cardinals and
either $\mu_N(\k)$ is inaccessible, or
$\theta < cf\mu_N(\k) < \k$.
\endproclaim

\proclaim{(10.1.4)  Remark}  If $\k$ is $\N\text{-controlled}$ but
$\k \not\in N$, then $\mu_N(\k)$ is singular.
\endproclaim

To see this, suppose that
$\mu > \k,\ \mu \in N$, with $\mu$ inaccessible.  Since $\k$ is $\N\text{-
controlled}$, let $\k \in d \in N$ where $d$ is an Easton set.
Thus, $sup\ d \cap \mu < \mu$ and clearly $sup\ d \cap \mu \in N$.
Now $\k \in d \cap \mu$, so $\k < sup\ d \cap \mu$.
But then,
it is easy to see that $\mu \neq \mu(\k)$, since if equality held, we would
have $N \cap [\k,\ \mu) = \emptyset$, contradicting that $\k < sup\ d \in
N \cap \mu$.
\bigskip
\noindent
{\bf (10.1.5)}
Now, let $\modseq$ be $(\Cal M,\ \th)\text{-standard}$ of length $\th + 1$.
For $i \leq \th$, let $\chi_i = \chi_{\Cal N_i},\ p\chi_i =
p\chi_{\Cal N_i}$.  Let $\Cal
N = \Cal N_\th = \bigcup\{\Cal N_i|i < \th\}$, and let $\chi = \chi_\th,
\ p\chi = p\chi_\th$, so
$dom\ \chi = \bigcup\{dom\ \chi_i|i < \th\}$, and for $\k \in dom\ \chi,
\ \chi(\k) = sup\ \{\chi_i(\k)|\k \in N_i\}$.  Also, for singular
cardinals, $\k \in [\th,\ \nu]$,
which are $\Cal N\text{-controlled},\ p\chi(\k) =
sup\ \{p\chi_i(\k)|i < \th\ \&\ \k\text{ is } N_i\text{-controlled}\}$.

Let $\k$ be a singular cardinal, $\k \in dom\ \chi$.
Note that since $cf\ \th = \th > \omega$, there is a club
$D \subseteq \th$ such that for all $i \in D,\ \chi_i(\k) \in C_{\chi(\k)}$.
This motivates the following.
\medskip
\proclaim{Definition}\ \ Let $\M,\ \th$ be as above,
and let $\modseq$ be $(\Cal M,\ \th)\text{-standard}$ of length $\th + 1$.
Let $\Cal N = \Cal N_\th$.
Let $N,\ N_i,\ \chi,\ p\chi,\ \chi_i,\ p\chi_i$ be as above.

Let $\k \geq \th$ be a singular cardinal, $\k \in N$.  $(\N_i|i \leq \th)$ is
$\bold {\M}\text{{\bf{-coherent} at }}\bold \k$
iff for all limit ordinals $\de \leq \th$
with $\k \in N_\de$, there is a club $D \subseteq \de$ such that
for all $i \in D,\ \chi_i(\k) \in C_{\chi_\de(\k)}$.
$(\N_i|i \leq \th)$ is
$\bold {\M}\text{{\bf{-coherent}}}$ if for all
singular cardinals $\k \in N \setminus \th,
\ (\N_i|i \leq \si)$ is $\M\text{-coherent at }\k$.
$(\N_i|i \leq \th)$ is {\bf strongly }
$\bold {\M}\text{{\bf{-coherent}}}$ iff for all
$i < \th$ and all singular cardinals $\k \in N_i,\ \chi_i(\k)
\in C_{\chi(\k)}$.  Finally,
$(\N_i|i \leq \th)$ is {\bf super }
$\bold {\M}\text{{\bf{-coherent}}}$ iff
$(\N_i|i \leq \th)$ is strongly $\M\text{-coherent}$ and for all
limit ordinals, $\si \leq \th$ and all
singular cardinals, $\k$ which are $\Cal N_\si\text{-controlled}$,
for sufficiently
large $i < \si,\ p\chi_i(\k) \in C_{p\chi_\si(\k)}$.
\endproclaim
\proclaim{(10.1.6)  Remark}
Let $(\N_i|i \leq \th)$ be
$(\Cal M,\ \th)\text{-standard}$ of length $\th + 1$.
For $i \leq \th$, let $\mu_i\ :=\ \mu_{\N_i}$, and let $\mu\ :=\ \mu_\th$.
Note that if $i < j$ then
$dom\ \mu_j \subseteq dom\ \mu_i$ and that if
$\k \in dom\ \mu_j$, then $\mu_j(\k)
\leq \mu_i(\k)$.  Thus, $dom\ \mu = \bigcap
\{dom\ \mu_i|i < \th\}$ and for
$\k \in dom\ \mu,\ \mu(\k)$ is the eventually
constant value of the $\mu_i(\k),\ i < \th$.
\endproclaim
\bigskip
\proclaim{(10.2) Lemma}  Let $\theta$ be regular, $\theta > \aleph_1$.  Let
$\nu  > cf\ \nu >> \theta$ be such that $(H_\nu,\in) \models$ a sufficiently
rich fragment of $ZFC$.  Let $\Cal M = (H_{\nu^+},\in,\cdots)$.
Let $C \subseteq
[H_{\nu^+}]^\theta$ be club.  Then there's strongly $\Cal M$-coherent
$(\Cal N_i|i \leq \theta)$, each $|\Cal N_i| \in C$.
\endproclaim

\demo{Proof} Without loss of generality, we may assume that
$X \in C \Rightarrow \Cal M|X \prec
\Cal M$.  We first build $(\Cal M_j|j \leq \theta^+)$,
an increasing continuous
tower of elementary submodels of $\Cal M$, each $|\Cal M_j| \in C,
\ \Cal M_j \in |\Cal M_{j+1}|,\ |\Cal M_{j+1}|$ closed under sequences
of length $< \theta$, for $j < \theta^+$.  Let $\chi_j = \chi_{\Cal M_j}$,\
$\chi = \chi_{\Cal M_{\theta^+}}$ be as in (10.1.5).

For singular $\kappa > \theta,\ \kappa \in |\Cal M_{\theta^+}|$, let
$E(\kappa) \subseteq \theta^+$ be club such that
$j \in E(\kappa) \Rightarrow \kappa
\in |\Cal M_j|$ and $\chi_j(\kappa) \in C_{\chi(\kappa)}$.  For $i <
\theta^+$, let \newline
$E_i = \bigcap\{E(\kappa)|\kappa \in |\Cal M_i|,\ \k > \th,
\ \k$ is singular$\}$, so
each $E_i$ is a club of $\theta^+$.  Let $E = \Delta_{i<\theta^+}
E_i =$ the diagonal intersection of the $E_i$.
Thus, \newline
$j \in E \Rightarrow (\forall
i < j)(\forall\kappa \in |\Cal M_i|)\chi_i(\kappa) \in
C_{\chi(\kappa)}$.  Let $E^\theta = \{\alpha \in E|cf\ \alpha = \th \}$ and
let $E^* = E^\theta \cup ((E^\th)^\prime \cap \theta^+)$ and let
$(j_i|i < \theta^+)$ be the increasing enumeration of $E^*$.  Thus for all
$i < \theta^+,\ cf\ j_{i + 1} = \theta$.  For $i \leq
\theta$, let $\Cal N_i = \Cal M_{j_i}$.  Then, $(\Cal N_i|i \leq \theta)$ is
strongly $\Cal M$-coherent.  All properties are clear from
construction, except possibly that for $i < \theta,
\ [|\Cal N_{i + 1}|]^{<\ \theta} \subseteq
|\Cal N_{i+1}|$.  This, however, is an easy consequence of the fact
that for successor $\zeta,\ |\Cal M_\zeta|$ is closed under sequences
of length $ < \theta$ and that $cf\ j_{i + 1} = \theta$.
\enddemo
\bigskip
\proclaim{Discussion}  We are now in a position to lay out the ideas behind
the proof, in (10.3) - (10.5), that if $(\Cal N_i|i \leq \th)$ is
strongly $\M\text{-coherent}$ then $(\Cal N_{\om i}|i \leq \th)$ is
super $\M\text{-coherent}$.
Let $\si \leq \th$ be a limit ordinal, and $\th < \k < \nu$ be a singular
cardinal.  We say that $\si$ is $\k\text{-good}$ if $\k$ is
$\N_\si\text{-controlled}$.  Now suppose that
$\de \leq \th$ is a limit of limit ordinals,
$\th < \k < \nu$ is a singular cardinal and that $\de$ is $\k\text{-good}$.
Let $\eta = p\chi_{\N_\de}(\k)$.
Our aim is to show that
for sufficiently large $\k\text{-good}$
limit ordinals $\si < \de,\ p\chi_{\N_\si}(\k) \in C_\eta$.

If we \lq\lq go up" to $\mu = \mu_{\N_\de}$ and let $\eta^\prime =
\chi_{\N_\de}(\mu)$, then, since $(\Cal N_i|i \leq \th)$ is
strongly $\M\text{-coherent}$, we have that for
$i < \de,\ \chi_{\N_i}(\mu) \in C_{eta^\prime}$.  Is there some way of
\lq\lq projecting" this fact back down to \lq\lq level $\k$"?
One such way would be to evaluate the $L\text{-scale}$ functions from
\lq\lq level $\mu$" (the $\hat f's$ at $\k$.
And, in fact, by (C) (1) of (1.2),
if we let $\eta^* = \hat f_{\eta^\prime}(\k),\
\{\hat f_\tau(\k)|\tau \in C_{\eta^\prime}\}$ will be a final
segment of $C_{\eta^*}$.  But what is the relationship between
$\eta$ and $\eta^*$, and, for $i < \de$, between
$p\chi_{\N_i}(\k)$ and $\hat f_{\chi_{\N_i}(\mu)}(\k)$.  The argument would
be complete, if we knew we had equality in the first case, and equality
in the second case for sufficiently large limit ordinals which are
$\k\text{-good}$.  This is exactly what will be proved in (10.4).
(10.3) supplies a technical result underlying the argument of (10.4).
In (10.5), we fill in the last few missing details of the above sketch,
in the presence of the result of (10.4).
\endproclaim
\bigskip
\noindent
{\bf (10.3)}
\medskip
If $\k$ is $\N\text{-controlled}$, set $g \in \Cal G_\N$ iff
$f \in |\N|,\ f$ is a function, $dom\ f$ is a set of $L\text{-cardinals}$
and for all $\tau \in dom\ f,\ \tau < f(\tau) < (\tau^+)^L$.  We also set
$\Cal G^L_\N\ :=\ \Cal G_\N \cap L$.

\proclaim{Proposition} If $\kappa$ is $\Cal N\text{-controlled},\
p\chi_{\Cal N}(\kappa) =
sup\ \{f(\kappa)|f \in \Cal G_\N,\ \k \in dom\ f\} =
sup\ \{f(\kappa)|f \in \Cal G^L_\N,\ \k \in dom\ f\}$.
\endproclaim

\demo{Proof} Clearly $p\chi_{\Cal N}(\kappa) \geq
sup\ \{f(\kappa)|f \in \Cal G_\N,\ \kappa \in dom\ f\} \geq
sup\ \{f(\kappa)|f \in \Cal G^L_\N,\ \k \in dom\ f\}$, so we show
that $p\chi_{\Cal N}(\kappa) \leq
sup\ \{f(\kappa) \in \Cal G_\N,\ \k \in dom\ f\}
\leq sup\ \{f(\kappa)|f \in \Cal G^L_\N,\ \k \in dom\ f\}$.
Since $\k$ is $\N\text{-controlled}$ (this is the whole point of the notion),
the last inequality is clear by covering,
so we prove the first.

Let $\xi < \kappa^+,\ \xi$ definable in $\Cal M$, by $\psi$, from
$x_1,\ \cdots,\ x_k \in |\Cal N|$ and $\kappa$.  Let $f(\tau) \simeq$
the least $\alpha < (\tau^+)^L$ such that $\Cal M \models
\psi(\alpha,\ x_1,\ \cdots,\ x_k,\ \tau)$, for
$L$-cardinals $\tau$.  Clearly
$f(\kappa) = \xi$ and for all $\eta \leq \nu,\ f|\eta \in |\Cal M|$.
Also, $\eta \mapsto f|\eta$ is $\Cal M$-definable.  Thus, if $\eta
\in |\Cal N|,\ f|\eta \in |\Cal N|$.  But clearly $\nu \in |\Cal
N|$.  Thus $f|\nu \in |\Cal N|$ and so $\xi =
(f|\nu)(\kappa)$.
\enddemo
\bigskip
\proclaim{Remarks}

(1)\ \ We could also have defined $\Cal G^{E}_\Cal N$ to be the set of
$f \in \Cal G_\Cal N$ such that $dom\ f$ is an Easton set, and
$\Cal G^{L,E}_\Cal N$ to be $\Cal G^E_\Cal N \cap L$, thereby
\lq\lq building in" the restriction to controlled $\k$.

(2)\ \ In connection with the alternative definition of the
$\chi_\Cal N$ and $p\chi_\Cal N$ as \lq\lq pressing down" functions, mentioned
at the end of (10.1.2), the above Proposition remains true, with these
alternative definitions, and the appropriately modified definition of the
various $\Cal G's$:  for $f \in \Cal G_\Cal N$ and $\k \in dom\ f,\
f(\k)$ would be required to be less than $\k$.
\endproclaim
\bigskip
\noindent
{\bf (10.4)\ \ }Suppose now that $\M$ is as in (10.2) and that
$\N^\prime \prec \M,\ card\ |\N^\prime| = \theta$ and let $\chi =
\chi_{\N^\prime},\ p\chi = p\chi_{\N^\prime}$
Let $\k$ be a
singular cardinal which is $\N^\prime\text{-controlled}$.
Let $\mu(\k) = \mu_{\Cal N^\prime}(\k)$, so that, by (10.1.4),
$\mu(\k)$ is a singular cardinal.
Let $\mu = \mu(\k)$,
let $\eta^\prime = \chi(\mu),\ \eta =
p\chi(\k)$ and suppose that $C_{\eta^\prime} \cap |\N^\prime|$ is cofinal in
$\eta^\prime$.  This will hold, in all cases of interest.


\proclaim{Lemma} $\eta = \hat f_{\eta^\prime}(\k)$.
\endproclaim

\demo{Proof}  We will end up applying (1.2)(B)(4)
(here, and in what follows, recall (9.3.2)!), so we must first show
that here, we have the hypothesis that
$(\hat f_\tau(\k)|\tau \in C_{\eta^\prime})$ is not eventually constant.
We begin with a number of easy observations, which we shall use
at various places in the proof.
\medskip
\roster
\item"{(1)}"  For $\tau$ a limit ordinal in $(\mu,\mu^+),
\ \hat f_\tau$ is canonically definable from $\tau$ in $\M$, so
for $\tau \in C_{\eta^\prime} \cap |\N^\prime|,\ \hat f_\tau \in |\N^\prime|$.
\item"{(2)}"  if $f,\ g \in |\N^\prime|$, where $f,\ g$ are functions with
domain the set of uncountable $L$-cardinals\ $<\ \mu$ and
$f <^* g$ then the least $\la_0$ such that
$(\forall\la \ge \la_0)f(\la) < g(\la)$ is definable from $f,\ g$ and is
therefore in $|\N^\prime|$.  So, since $\la_0 \in
|\N^\prime|$ and $\la_0 < \mu$, we must have  $\la_0 < \kappa$.
\endroster
\medskip

We are now in a position to argue that
$(\hat f_\tau(\k)|\tau \in C_{\eta^\prime})$ is not eventually constant.
We will
do this by proving that for $\tau_1 < \tau_2$, both in $C_{\eta^\prime} \cap
|\Cal N^\prime|$, for all $\la \geq \k,\ \hat f_{\tau_1}(\la)
< \hat f_{\tau_2}(\la)$.
In particular, this means that the map from
$C_{\eta^\prime} \cap |\Cal N^\prime|$ to $\Phi$, given by $\tau \mapsto
\hat f_\tau(\k)$ is order preserving, so $\Phi$ has limit order type, as
required, since clearly $C_{\eta^\prime}$ does.  So, suppose $\tau_1,\ \tau_2$
are as above.  Applying (1), we have that for $i = 1,\ 2,\ \hat f_{\tau_i} \in
|\Cal N^\prime|$.  But then, we have the desired conclusion, by applying (2),
with $f = \hat f_{\tau_1},\ g = \hat f_{\tau_2}$.

As we have just proven, we have the hypotheses of (1.2)(B)(4),
so, by (1.2)(B)(4), $\hat f_{\eta^\prime}(\k) =
sup\ \{\hat f_\tau(\k)|\tau \in
C_{\eta^\prime}\}$.  Further,
by (1.2)(B)(2),
$sup\ \{\hat f_\tau(\k)|\tau \in
C_{\eta^\prime}\}  = sup\ \{\hat f(\k)|\tau \in C_{\eta^\prime} \cap
|\Cal N^\prime|\}$.
Again, by (1), if $\tau$ is as in (1),\
$\hat f_\tau(\k) < p\chi(\kappa)$
so finally, $\hat f_{\eta^\prime}(\k) \leq \eta$.

Clearly $\eta = sup\ \{f(\k)|f \in |\N^\prime| \cap L_{\mu^+}\} =
sup\ \{f(\k)|f \in |\N^\prime| \cap L_{\eta^\prime}\}$.  Thus, it
suffices to show:
\medskip
\roster

\item"{($\ast$):}" if $f \in |\N^\prime| \cap L_{\eta^\prime}$,
$dom\ f$ is the set of uncountable
$L$-cardinals $<\ \mu$ and for $\la \in dom\ f,
\ f(\la) \in (\la,(\la^+)^L)$,
then there's $\gamma \in |\N^\prime| \cap L_{\eta^\prime}$ such that
\itemitem {(a)} $f <^* \hat f_\gamma$, and
\itemitem {(b)} for all $\lambda \geq \k,\ f(\la) < \hat f_\gamma(\la)$.
\endroster
\medskip
\noindent
Now, the existence of a
$\gamma \in |\N^\prime| \cap L_{\eta^\prime}$ satisfying
(a) is an easy consequence of $\N^\prime \prec \M$ and the fact, which holds in
$\M$, that $(\hat f_\xi|\xi \in (\mu,\mu^+))$ is an $L\text{-scale}$, by
(B)(5) of (1.2).  But
then for such a $\gamma,\ \hat f_\gamma \in |\N^\prime|$,
and then (b) follows immediately from (2), with $g = \hat f_\gamma$.
\enddemo
\bigskip

\proclaim{(10.5)  Lemma} If $(\N_i|i \leq \th)$ is strongly $\M$-coherent then
$(\N_{\om i}|i \leq \th)$ is super $\M$-coherent. \endproclaim

\demo{Proof}  We fill in the details of the argument sketched in
the Discussion following (10.2).  We adopt the notation and terminology
established there.  Let $\de \leq \th$ be a limit of ordinals.  Suppose
that $\th < \k < \nu$ and that $\de$ is $\k\text{-good}$.  By (10.1.6),
there is $i_0 < \de$ such that if $i_0 \leq i < \de,\ \mu_{\N_i}(\k) =
\mu_{\N_\de}(\k)$.  Let $\mu, \eta,\ \eta^\prime,\ \eta^*$
be as in the Discussion.

Now, let $\N^\prime = \N_\de$.  Since $(\N_i|i \leq \th)$ is
strongly $\M\text{-coherent}$, it is easy to see that
$C_{\eta^\prime} \cap |\N^\prime|$ is cofinal in $\eta^\prime$, so we have
the hypotheses of (10.4).  Thus,
by the Lemma of (10.4),
$\eta = \hat f_{\eta^\prime}(\k)\ (\ =\ \eta^*)$.

Suppose, now that $i_0 < \si < \de$, where $\si$ is a $\k\text{-good}$
limit ordinal.  Since $i_0 < \si,\ \mu_{\N_\si}(\k) = \mu$.
Therefore, we can apply (10.4), again,
but with $\N^\prime = \N_\si$;
just as in the preceding paragraph, this give us that
$p\chi_{\N_\si}(\k) = \hat f_{\chi_{\N_\si}(\mu)}(\k)$.
The conclusion is now clear, as in the Discussion:
$\{\hat f_\tau(\k)|\tau \in C_{\eta^\prime}\}$ is a final segment
of $C_\eta$, for all $i < \de,\ \chi_{\N_i}(\mu)
\in C_{\eta^\prime}$, and for all
$\k\text{-good}$ limit ordinals, $\si$, with
$i_0 < si < \de,\ p\chi_{\N_\si}(\k) = \hat f_{\chi_{\N_\si}(\mu)}(\k)$,
so for all sufficiently large $\k\text{-good}$ limit ordinals, $\si$,
with $i_0 < \si < \de,\ p\chi_{\N_\si}(\k) \in C_\eta$, as required.
\enddemo
\bigskip
Now, clearly, combining (10.2) and (10.5), we have proved Lemma 3.

\bigskip
\noindent
{\bf (10.6)}   We now expand somewhat on the proof of (10.4).  We have
already noted that $\mu(\k)$ is a limit of limit cardinals.  Suppose first
that $\k$ is of the form $\la^{+\om}$.  Then, for all such $\la$, and all
$\N^\prime \prec \M,\ \k \notin
|\N^\prime| \Rightarrow [\la,\ \k] \cap |\N^\prime| = \emptyset$.
Thus, in this setting, for all $\la \in D_\k,\ \mu(\la) = \mu(\k)$.

If $\k$ is a singular limit of limit cardinals and
$\k \cap |\N^\prime|$ is bounded
in $\k$, then, on a tail of $D_\k,\ \mu(\la) = \mu(\k)$.  Let us then examine
the most difficult case, where $\k$ is a singular limit of limit cardinals,
$\k \not\in |\N^\prime|$, but
$\k \cap |\N^\prime|$ is cofinal in $\k$; note that if
$|\N^\prime|$ is closed for sequences of length
$< \theta$ (as was the case, in
the context of (10.5), taking $\N^\prime = \N_\th$),
this means that $cf\ \k =
\th$.  Note also that we may even have $D_\k \subseteq |\N^\prime|$.
Recall that, in this
latter case, $p\chi(\la) = \chi(\la)$, for $\la \in D_\k$.  Even if $\la \in
D_\k \setminus |\N^\prime|$, we still have $\mu(\la) < \k < \mu(\k)$.

Our principal aim is to show that one inequality of the Lemma of (10.4) remains
true when we replace $\eta = p\chi(\k)$ by $\sigma = p\chi(\la)$ and
$\hat f_{\eta^\prime}(\k)$ by $\hat f_{\eta^\prime}(\la)$,
but maintain $\eta^\prime =
\chi(\mu(\k))$, instead of using $\sigma^\prime = \chi(\mu(\la))$.  Of course,
by (10.4) with $\la$ in place of $\k$, we \underbar{do} have $\sigma =
\hat f_{\sigma^\prime}(\la)$, and this is our point of departure in proving:

\proclaim{Lemma} $\sigma \geq \hat f_{\eta^\prime}(\la)$,
on a tail of $\la \in D_\k$.
\endproclaim

\demo{Proof}  We follow the proof of (10.4).
Obtaining $\leq$ seems problematical since the
proof of the analogue of ($\ast$) of (10.4) does not seem to go over.

First, take $\la$ sufficiently large that
$\mu^*_1(\eta^\prime) < \la$, where $\mu^*_1(\eta^\prime)$ is as in (5.3).
This is
possible, since, in (10.4), we showed that
$\mu^*_1(\eta^\prime) < \k$.  Now, for
such $\la$, the proof in (10.4), that
$\hat f_{\eta^\prime}(\k) \leq \eta$, goes over
verbatim to show that $\hat f_{\eta^\prime}(\la) \le \sigma$. \enddemo
\bigskip
\subheading{\S 11.  ODDS AND ENDS}
\medskip
We close by providing the proofs of two small results needed for
\cite{5}.  In (11.1) we prove the Proposition of (1.5) of
\cite{5} needed for the construction of the very tidy conditions.
In (11.2) we show, as promised in (2.1.1) of
\cite{5}, that,
without loss of generality, the system of
$b_\a$ for $\a$ which are
multiples of $card\ \a$, which is inaccessible,
can be taken to be {\it tree-like}.
\bigskip
\proclaim{(11.1)  Proposition}
Let $\theta > \aleph_1$ and let $\nu,\ \Cal M$ be as
in (10.2).  Let $d \subseteq [\theta,\nu)$ be an Easton set of cardinals, and
let $\ga$ be a function with domain $d$ such that for all $\k \in d,\
\ga(\k) < \k^+$.  Then, there is a function $\ga^*$ with domain $d$ such that
for all $\k \in d$ which are either singular or of the form
$\aleph_\tau$, with $\tau > 1$ and odd,
$\ga^*(\k) > \ga(\k)$ and such that for all
singular $\k \in d$, letting $\a = \ga^*(\k),\ f^*_\a =^*\ \ga^*|D_\k$.
Further, if $\Cal N \prec \Cal M$ with $(\theta + 1) \cup \{\ga\} \subseteq
|\Cal N|$, then $\ga^* \in |\Cal N|$.
\endproclaim
\demo{Proof}
We first define a function $\ga_1$ as follows:
if $\th \leq \k < \nu$, where $\k$ is of the
form $\aleph_{\a + \omega}$, we let
$\ga_1(\k) =\ $the least $\eta \in (\ga(\k),\
\k^+)$ such that $\eta$ is a multiple of $\k^2$ and
$\ga|[\aleph_\a,\ \k) <^*
\hat f_\eta$ (where we define $\ga(\la)$ to be the usual default
value, $\la^22$, for $\la \in [\aleph_\a,\ \k) \setminus d$).
For regular cardinals $\la \in d \cap [\th,\ \nu)$, we let $\ga_1(\la) =
max(\ga(\la) + \la^2,\ \hat f_{\ga_1(\la^{+\omega})}(\la))$.
For all other $\k \in
[\th,\ \nu)$, we let $\ga_1(\k) = \ga(\k) + \k^2$.\par

To obtain $\ga^*$ from $\ga_1$,
we first define, by recursion on $n < \om$, ordinals, $\nu_n,\ \eta_n$,
and a function, $f_n$.  We will have that if $\nu_n > \th$, then
$\nu_{n+1} < \nu_n$, so there will be $m < \om$ such that
$\nu_{m+1} \leq \th < \nu_m$.  We stop the recursion at this $m$.

Let $\nu_0 = \nu$.  Having defined $\nu_n$,
if $\nu_n$ is a singular limit of limit cardinals with $\nu_n \in (\th,\
\nu]$, we let $\eta_n \in [\ga_1(\nu_n),\ \nu^+_n)$ be the least
$\eta$ which is a multiple of $\nu^2_n$ such that:
$$(\ast)_n:  \hat f_\eta >^* \ga_1|d \cap \nu_n$$
and we let $f_n = \hat f_{\eta_n} \cup \{(\nu_n,\
\eta_n)\}$.
Once again, this is possible by Covering, because we have
taken the precaution of restricting to an Easton, $d$.

So, having defined $\nu_n,\ \eta_n,\ f_n$, satisfying
$(\ast)_n$, we
define: \par
{\narrower\smallskip\noindent  $\nu^0_{n + 1} =\ $
the least cardinal $\nu^\prime \in [\th,\ \nu_n)$ such
that for singular $\k \in d \cap [\nu^\prime,\ \nu_n),\ \ga_1(\k) <
f_n(\k)$. \medskip }

\noindent
Having defined $\nu^i_{n + 1}$, if $\nu^i_{n + 1} \leq \th$, we set
$\nu_{n + 1} = \th,\ m = n$ and we stop.  If
$\nu^i_{n + 1} > \th$ is a singular limit of
limit cardinals, we set $\nu_{n + 1} = \nu^i_{n + 1}$.
If $\nu^i_{n + 1} > \th$ is either a successor cardinal or of the form
$\aleph_{\tau + \omega}$, we set $\nu^{i + 1}_{n + 1} =\ $ the largest
limit of limit cardinals\ $<\ \nu^i_{n + 1}$.  Finally, if $\nu^i_{n + 1} >
\th$ is inaccessible, we set
$\nu^{i + 1}_{n + 1} = sup\ d \cap \nu^i_{n + 1}$.

Clearly there is $i < \omega$ such that $\nu_{n + 1} = \nu^i_{n + 1}$ and
either $\nu^i_{n + 1} \leq \th$ or $\nu^i_{n + 1}$ is a singular limit of
limit cardinals.  In all cases, we let
$a_n = \{\k \in (\nu_{n + 1},\ \nu^0_{n + 1}]|\k
\ \text{is a singular cardinal} \}$; note that $a_n \cap d$ is finite,
and for all $\k \in a_n \cap d,\ \k$ is not a limit of limit cardinals.
When $\nu_{n + 1} > \th$, we have $m > n$, and we continue, to define
$\eta_{n + 1}$ and $\nu_{n + 2}$.  Clearly $m < \omega$, i.e., for some
$n,\ \nu_{n + 1} \leq \th$.

We now define
$\ga^*$: \par
{ \narrower\smallskip\noindent
if $\k \in d,\ \k$ is singular, $\k \not\in\
\bigcup\{a_n|n \leq m\}$,
we set $\ga^*(\k) =  f_n(\k)$, where $n$ is such that $\nu_{n
+ 1} < \k \leq \nu_n$.  If $\k \in a_n$, where $n \leq m$,
we let $\ga^*(\k) = \ga_1(\k)$.  Finally, if $\la = \aleph_\tau$,with  $\tau$
odd, $\la \in d$, we set
$\k = \aleph_{\tau + \omega}$ and we set
$\ga^*(\la) = max(\ga_1(\la),\ \hat f_{\ga^*(\k)}(\la))$.  For all other
$\la \in d$, we set $\ga^*(\la) = \ga_1(\la)$.
\medskip }

\noindent
It is clear that:

{\nar (\#)  If $\M$ is as in (10.1), $\N^\pr \prec \M,\ N^\pr\ :=\ |\N^\pr|,
\ \th + 1 \subseteq N^\pr,\ card\ N^\pr = \th,
\ [N^\pr]^{<\ \th} \subseteq N^\pr,
\ d,\ \ga \in N^\pr$, then $\ga^* \in N^\pr$.
\medskip }

But then, clearly, $\ga^*$ is as required.
\enddemo
\bigskip
\noindent
{\bf (11.2)  GETTING \lq\lq TREE-LIKE" } $\bold{b_\a.}$
\medskip
We begin by recalling some notions from the Introduction and (2.1) of
\cite{5}.  First, recall that for inaccessible $\k,\ U(\k)$ is the set of
multiples of $\k$ in $(\k,\ \k^+)$.
Let $\k$ be inaccessible.  Recall that a system,
$(b_\a|\a \in U(\k))$ of almost-disjoint cofinal subsets of $\k$
was called {\it decodable} if

$$
\split
(\ast):\ \ \ \text{  for all } \theta\in (\kappa,\ \kappa^+),
\ (b_\alpha|\alpha \leq \theta) \in L[A \cap \theta],\\
\text{ and is \lq\lq canonically definable" there}.
\endsplit
$$

Recall that Corollary 4 of the Introduction of \cite{5} gives that
for all inaccessible $\k$,
there is decodable
$\overset{\rightharpoonup}\to{b} =
(b_\alpha|\alpha\in U(\k))$ of cofinal
almost-disjoint subsets of $\kappa$
as above.

In (2.1.1) of \cite{5}, we defined
$\Cal U\ :=\ \bigcup\{U(\k)|\k \text{ is inaccessible}\}$,
and we considered the following additional property of the
system $(b_\eta|\eta \in \Cal U)$ which we called
{\bf tree-like }:
\medskip

{\nar whenever
$\eta_1,\ \eta_2 \in \Cal U$, if $\xi \in b_{\eta_1} \cap b_{\eta_2}$,
then $b_{\eta_1} \cap \xi = b_{\eta_2} \cap \xi$.\medskip}

We promised there, to show, here:
\proclaim{Lemma}
Without loss of generality, we can
assume that $(b_\eta|\eta \in \Cal U)$ is tree-like and has the following
additional property:  $b_\eta = range\ g_\eta$, where $g_\eta$ is a function,
$dom\ g_\eta = \{\aleph_\tau|\aleph_\tau < card\ \eta\ \&\
\tau \text{ is an even successor ordinal}
\}$; further, for all $\xi \in b_\eta,\ \xi$ is  a multiple of 4 but
not of 8.
\endproclaim
\demo{Proof}
This is actually a rather simple observation; for the record, the following is
one way this can be achieved.

For inaccessible $\kappa$ and $\alpha \in U(\kappa)$, and
$\lambda < \kappa$ of the form $\aleph_\tau$, where
$\tau$ is an even successor, let
$\zeta_\alpha(\lambda)$ be the rank of $b_\alpha \cap \lambda$ in
$<_{L[A \cap \lambda^+]}$, and let $g_\alpha(\lambda) = $
the $\zeta_\alpha(\lambda)^{th}\ \eta$ such that
$\lambda < \eta < \lambda^+$ and $\eta$ is a multiple of
4 but not of 8.  Then let $b^*_\alpha = range\ g_\alpha$.
It is clear that the $b^*_\alpha$ are decodable, since
the $b_\alpha$ were, and that they have the desired tree-like property.
\enddemo

\enddocument